\documentclass[preprint,12pt]{elsarticle}

\usepackage{amssymb}

\usepackage[utf8]{inputenc}
\usepackage[T1]{fontenc}
\usepackage{hyphsubst}
\usepackage[margin=2cm]{geometry}
\usepackage{lmodern}
\usepackage[english]{babel}
\usepackage{mathtools, nccmath}

\usepackage{blindtext} 
\usepackage{charter} 
\usepackage{graphicx, multicol} 
\usepackage{xcolor} 
\usepackage{pseudocode} 
\usepackage{parskip}
\usepackage{subfig}
\allowdisplaybreaks

\usepackage{bef_aras_nothm,todonotes}
\usepackage{amsmath,amsthm,amssymb}
\usepackage{fixmath}
\usepackage{upgreek}
\usepackage{amsfonts}
\usepackage{microtype}
\usepackage{float}
\usepackage{mathrsfs}
\usepackage{hyperref}
\usepackage{paralist}
\hypersetup{
     colorlinks   = true,
     citecolor    = blue
}
\usepackage[nameinlink]{cleveref} 
\crefname{figure}{\figurename}{\figurename}  
\usepackage{tikz}
\usetikzlibrary{shapes,arrows,positioning, calc}
\usepackage{verbatim}
\usepackage{afterpage}
\usepackage{glossaries}

\tikzset{
>=stealth',
help lines/.style={dashed, thick},
axis/.style={<->},
important line/.style={thick},
connection/.style={thick, dotted},
}

\theoremstyle{plain}
\newtheorem{theorem}{Theorem}[section]
\newtheorem{lemma}[theorem]{Lemma}
\newtheorem{proposition}[theorem]{Proposition}
\newtheorem{definition}[theorem]{Definition}
\newtheorem{corollary}[theorem]{Corollary}
\theoremstyle{definition}

\newtheorem{remark}[theorem]{Remark}

\newcounter{ArasCounter}

\newcounter{BeatriceCounter}

\makeatletter
\g@addto@macro{\thm@space@setup}{\thm@headpunct{}}
\makeatother

\numberwithin{equation}{section}

\journal{Journal of Computational Physics}

\begin{document}

\begin{frontmatter}



\title{Error Estimation for Physics-informed Neural Networks Approximating Semilinear Wave Equations}


\author[inst1]{Beatrice Lorenz}

\affiliation[inst1]{organization={Department of Mathematics, Ludwig Maximilians Universität München},
            postcode={80798}, 
            country={Germany}}

\author[inst3]{Aras Bacho\corref{cor1}}
\author[inst1,inst2]{Gitta Kutyniok}

\affiliation[inst2]{organization={Munich Center for Machine Learning (MCML)},
            city={Munich},
            country={Germany}}
            
\affiliation[inst3]{organization= {Department of Computing and Mathematical Sciences, California Institute of Technology},
            city={Pasadena},
            country={CA, USA}}

\cortext[cor1]{Corresponding author, E-mail: bacho@caltech.edu}

\begin{abstract}
This paper provides rigorous error bounds for physics-informed neural networks approximating the semilinear wave equation. We provide bounds for the generalization and training error in terms of the width of the networks layers and the number of training points for a tanh neural network with two hidden layers. Our main result is an error bound for the $H^1([0,T];L^2(\Omega))$-norm of the error in terms of the training error and the number of training points, which can be made arbitrarily small under some assumptions. This norm also includes a bound on the total error. We illustrate our theoretical bounds with numerical experiments.
\end{abstract}



\end{frontmatter}


\section{Introduction}
Partial differential equations (PDEs) play an essential role in describing numerous physical phenomena across various scientific disciplines, ranging from fluid dynamics and electromagnetism to quantum mechanics and finance. Solving these equations analytically is often challenging or even impossible, necessitating the utilization of other methods to obtain approximate solutions.

One way to find approximate solutions to partial differential equations is via classical numerical methods. These methods have been studied for years and already have strong theoretical foundations when it comes to error estimation \cite{Golub}. However, in recent years with the rise of machine learning as a whole, there has also been an increased interest in applying machine learning methods to the problem of finding approximate solutions to PDEs. Among these methods are physics-informed neural networks (PINNs) which were conceived as feed-forward neural networks that incorporate the dynamics of the PDE into their loss function \cite{Raissi}. 

In this article, we will be providing error bounds for physics-informed neural networks approximating semilinear wave equations. The semi-linear wave equation describes the propagation of waves in various physical systems and can account for nonlinear effects such as self-interaction, interaction with an external force, energy dissipation and amplitude reduction. It has a wide range of applications in physics, engineering, and other fields \cite{Galstian,Conolly,Struwe}. PINN approximations of linear wave equations have previously been studied in \cite{Mosely,WangWaveEq}. These results show that that PINNs are able to solve the wave equation with high accuracy. Nevertheless, it is noted that the wave equation presents unique challenges for computational solvers due to the multi-scale, propagating and oscillatory nature of its solutions. So far, there are no results for PINNs approximating the semilinear wave equations and no error estimates for physics-informed neural networks approximating wave equations in general.

To remedy this, we will adapt an approach introduced by De Ryck et. al. \cite{DeRyck1} that has already been used successfully to find error bounds for PINNs approximating the Navier-Stokes \cite{DeRyck1} and the Kolmogorov equation \cite{DeRyckKolmogorov}. The semilinear wave equation presents two additional challenges that none of the previously mentioned equations posses: firstly, it includes the second time-derivative and secondly, it involves a general non-linearity which makes obtaining error bounds more challenging. Finding a good error bound is not only an important goal in itself, providing a valuation for the quality of the PINN approximation, but also serves as an appraisal of the methodology proposed in \cite{DeRyck1}, possibly warranting further application to other partial differential equation.

\subsection{PINNs Explained}\label{PINNsExplained}
Physics-informed neural networks are (deep) neural networks that approximate PDE solutions by minimizing a loss function that includes terms reflecting the initial and boundary conditions, as well as the PDE dynamics at selected points in the domain. The vanilla PINN, first proposed by Raissi et. al. \cite{Raissi} is a variation of a classic feed forward neural network with an amended loss function. 



To introduce the basic idea behind physics-informed neural networks, we consider the following abstract PDE framework:
\begin{equation}\label{DeRyckBVP}
    \begin{split}
        &\mathcal{D}[u](x,t)=0,\\
        &\mathcal{B}[u](y,t)=\psi(y,t),\\
        &u(x,0)=\varphi(x),\qquad\text{for }x\in\Omega,y\in\partial\Omega,t\in[0,T],
    \end{split}
\end{equation}
where $\Omega\subset\mathbb{R}^d,d\in\mathbb{N}$ is compact, $T>0$ and $\mathcal{D}, \mathcal{B}$ are the differential and boundary operators, $u:\Omega\times[0,T]\rightarrow\mathbb{R}^m,m\in\mathbb{N}$ is the solution of the PDE, $\psi:\partial\Omega\times[0,T]\rightarrow\mathbb{R}^m$ specifies the (spatial) boundary conditions and $\varphi:\Omega\rightarrow\mathbb{R}^m$ are the initial conditions.

The \textit{pointwise residuals} for any sufficiently smooth function $f:\Omega\times[0,T]\rightarrow\mathbb{R}^m$ are defined as as
\begin{equation}\label{Residuals}
    \begin{split}
        &\mathcal{R}_i[f](x,t)=\mathcal{D}[f](x,t),\\
        &\mathcal{R}_s[f](y,t)=\mathcal{B}[f](y,t)-\psi(y,t),\\
        &\mathcal{R}_t[f](x)=f(x,0)-\varphi(x),\qquad x\in\Omega,y\in\partial\Omega,t\in[0,T].
    \end{split}
\end{equation}
The residuals measure how well a function satisfies the boundary value problem \ref{DeRyckBVP}. For the exact solution, all residuals equal zero. The PINN framework seeks to minimize the difference between a neural network solution $u_{\theta}$ and the actual solutions $u$. This difference can be expressed in terms of these residuals and is called the \textit{generalization error} in \cite{DeRyck1}:
\begin{equation}\label{Generalization error}
    \mathcal{E}_G(\theta)^2=\int_{\Omega\times[0,T]}|\mathcal{R}_i[u_{\theta}](x,t)|^2dxdt+\int_{\partial\Omega\times[0,T]}|\mathcal{R}_s[u_{\theta}](x,t)|^2ds(x)dt+\int_{\Omega}|\mathcal{R}_t[u_{\theta}](x)|^2dx.
\end{equation}
In practice these integrals are usually evaluated using a suitable numerical quadrature rule leading to the definition of the \textit{training error} $\mathcal{E}_T(\theta,\mathcal{S})$:
\begin{equation}\label{TrainingErrors}
    \begin{split}
        \mathcal{E}_T(\theta,\mathcal{S})^2&=\mathcal{E}_T^i(\theta,\mathcal{S}_i)^2+\mathcal{E}_T^s(\theta,\mathcal{S}_s)^2+\mathcal{E}_T^t(\theta,\mathcal{S}_t)^2\\
        &=\sum_{n=1}^{N_i}w_i^n|\mathcal{R}_i[u_{\theta}](x_i^n,t_i^n)|^2+\sum_{n=1}^{N_s}w_s^n|\mathcal{R}_s[u_{\theta}](x_s^n,t_s^n)|^2+\sum_{n=1}^{N_t}w_t^n|\mathcal{R}_t[u_{\theta}](x_t^n)|^2,\\  
    \end{split}
\end{equation}
with quadrature points $\mathcal{S}_i=\{(x_i^n,t_i^n)\}_{n=1}^{N_i}, N_i\in\mathbb{N}, \mathcal{S}_s=\{(x_s^n,t_s^n)\}_{n=1}^{N_s}, N_s\in\mathbb{N}, \mathcal{S}_t=\{(x_t^n)\}_{n=1}^{N_t},N_t\in\mathbb{N}$ and the corresponding quadrature weights $w_i^n, w_s^n, w_t^n\in\mathbb{R}_+$.

Depending on the PDE in question the evaluation of the loss terms involves computing the gradients and higher derivatives of $u_{\theta}$, which is usually done using automatic differentiation. Furthermore, due to the computation of second order derivatives it is crucial to choose the activation function in a PINN framework carefully. Smooth activation functions, such as the sigmoid and hyperbolic tangent, can be used to guarantee the regularity of PINNs, allowing estimates of PINN generalization error to be accurate \cite{Mishra1}. 

Besides these so called vanilla PINNs, numerous variants have been developed over the years, each attempting to address specific problems within the original PINN framework or optimizing it for specific problems, such as VPINNs \cite{Kharazmi}, CAN-PINNs \cite{Chiu}, XPINNs \cite{Jagtap2}, cPINNs \cite{Jagtap1} and more.  Outside of PINNs, another popular approach to solving PDEs with machine learning are Kernel-based methods \cite{chen2021solving,bacho26} and so-called neural operator networks which aim to learn the underlying solution operator that maps initial and boundary conditions to the solution for a whole class of PDEs. This is contrasted with PINNs where the goal is to find the solution to a single PDE. Operator networks include the Fourier Neural Operator \cite{li2}, DeepONets \cite{Lu2} and other architectures \cite{Li2022FourierNO,Li2021PhysicsInformedNO,Li2020NeuralOG,Tripura2022WaveletNO, Viswanath, Raonic2023ConvolutionalNO, bach25}.

\subsection{Error Estimation Results for PINNs}

One of the first results concerning the convergence and error estimation of physics-informed neural networks is due to Shin et. al. \cite{Shin1,Shin2}. In \cite{Shin1} they show a probabilistic bound of the total error in terms of the training error as well as the strong convergence of the PINN solution to the ground truth as the number of collocation points goes to infinity when using a H\"older regularized loss function. In \cite{Shin1} they do this for a set of generic linear second order elliptic and parabolic type PDEs and in \cite{Shin2} a similar result, including more quantitative estimates, is provided for advection equations. Additionally, the results by  Mishra and Molinaro \cite{Mishra2} give an abstract framework for a bound of the total error in terms of the training error, the number of training points and stability bounds for a general abstract class of partial differential equations. They apply their framework to the concrete cases of viscious scalar conservation laws and incompressible Euler equations \cite{Mishra2}.

In \cite{DeRyck1} De Ryck et. al. expand the framework of \cite{Mishra2} by also including bounds on the training and generalization error and show the existence of two-layer tanh neural networks with small generalization and training errors. They introduce a general methodology for estimating the total error incurred by a physics-informed neural network that can be transferred to other partial differential equations. In \cite{DeRyck1} they provide error bounds for the Navier-Stokes equation and in \cite{DeRyckKolmogorov} they use a similar approach to find error bounds for PINNs approximating a class of Kolmogorov PDEs. Their approach was also successfully used for the primitive equations in \cite{Hu}.

\subsection{Methodology}
A key goal of the PINN approach is to find a neural network solution $u_{\theta}$ to \ref{DeRyckBVP} such that the total error between the computed PINN solution and the actual solution $\|u_{\theta}-u\|$ is small. In this paper, we will adress the following three questions:
\begin{itemize}\label{3Q}
    \item[Q1.] Given a tolerance $\varepsilon>0$, do there exist neural networks $u_{\hat{\theta}}, u_{\tilde{\theta}}$, parametrized by $\hat{\theta}, \tilde{\theta}$ such that the corresponding generalization and training errors are small, i.e. $\mathcal{E}_G(\hat{\theta}),\mathcal{E}_T(\tilde{\theta},\mathcal{S})<\varepsilon$?
    \item[Q2.] Given a PINN $u_{\theta}$ with a small generalization error, is the corresponding total error $\|u-u_{\theta}\|$ small, i.e.$\|u-u_{\theta}\|<\delta(\varepsilon)$, for some $\delta(\varepsilon)\sim\mathcal{O}(\varepsilon)$?
    \item[Q3.] Given a small training error $\mathcal{E}_T(\theta^*,\mathcal{S})$ and a sufficiently large training set $\mathcal{S}$, is the corresponding generalization error $\mathcal{E}_G(\theta^*)$ also proportionally small?
\end{itemize}
Q1 shows that in principle, the PDE residuals as well as the training and generalization error of a physics-informed neural network can be made arbitrarily small, since there exist neural networks where exactly that is the case. An affirmative answer to Q2 ensures that, at least theoretically, we can find a PINN solution such that the total error is arbitrarily small when the generalization error is small.Consequently, an affirmative answer to both Q1 and Q2 implies that in principle, there is a neural network such that the total error is appropriately small. The question that remains to be answered is whether we can find such a network in practice. That is where question three comes in. An affirmative answer to Q3 entails that if we get a small training error in practice, the generalization error will also be small and together with Q2 this implies a small total error. Together with an answer to Q1 which entails that there exists a PINN with an arbitrarily small training error, this leads to the guaranteed existence of a PINN with what can be assumed is a reasonably small total error. Moreover, under the assumption that the PINN finds the global minimum, network architecture and training set specifications can be given such that the total error of the network can be made arbitrarily small, i.e. the total error will be $\mathcal{O}(\varepsilon)$ for any $\varepsilon>0$.

\subsection{Contributions}

In the following we summarize our contributions:
\begin{itemize}
    \item We show that there exists a neural networks such that their generalization and training errors are arbitrarily small, thereby answering Q1. Moreover, we give explicit bounds to the number of layers and neurons in each layer to achieve this goal.
    \item We affirmatively answer Q2 by proving that the sum of the total error and the total error of the time derivative can be bound in terms of the PDE residuals.
    \item We show that the sum of the total error and the total error of the time derivative can be bound in terms of the training error and the size of the training set. Moreover, under the assumption that the algorithm finds the global minimum we give specifications for the network architecture and the size of the training set such that the total error, plus the total error of the time derivative can be made arbitrarily small.
    \item To sum up, we have found an \textit{a priori} as well as \textit{a posteriori} bound for the sum of the the total error and the total error of the time derivatives between the exact solution of a semilinear wave equation and the PINN solution.
    \item We have improved on previous results in this area by not only bounding the total error, but the total error plus the total error of the time derivative, ensuring a better approximation.
    \item We have validated these bounds for a damped wave equation with numerical experiments.
    \item We have dealt with the additional difficulties of the second time derivative in the semilinear wave equation as well as rather general nonlinearities $f$.
\end{itemize}

\section{Mathematical Preliminaries}

The following definition of a neural network will be used throughout:
\begin{definition}[Neural Network \cite{DeRyck1}]\label{NeuralNetwork}
    Let $R\in(0,\infty], L,W\in\mathbb{N}$ and $l_0,...,l_L\in\mathbb{N}.$ Let $\sigma:\mathbb{R}\rightarrow\mathbb{R}$ be a twice differentiable activation function and define
    \begin{equation*}
        \Theta=\Theta_{L,W,R}:=\bigcup_{L'\in\mathbb{N},L'\leq L}\bigcup_{l_0,...,l_L\in\{1,...,W\}}\bigtimes_{k=1}^{L'} \left( [-R,R]^{l_k\times l_{k-1}}\times[-R,R]^{l_k} \right)
    \end{equation*}
    For $\theta\in\Theta_{L,W,R}$ we define $\theta_k:=(W_k,b_k)$ and $\mathcal{A}_k:\mathbb{R}^
    {l_{k-1}}\rightarrow\mathbb{R}^{l_k}:x\mapsto W_kx+b_k$ for $1\leq K\leq L$ and we define $f^{\theta}_k:\mathbb{R}^
    {l_{k-1}}\rightarrow\mathbb{R}^{l_k}$ by
    \begin{equation*}
        f^{\theta}_k(z)=
        \begin{cases}
            \mathcal{A}^{\theta}_L(z) & k=L,\\
            (\sigma\circ\mathcal{A}^{\theta}_k)(z) & 1\leq k<L.
        \end{cases}
    \end{equation*}
    We denote by $u_{\theta}:\mathbb{R}^{l_0}\rightarrow\mathbb{R}^{l_L}$ the function that satisfies for all $z\in\mathbb{R}^{l_0}$ that
    \begin{equation*}
        u_{\theta}(z)=(f_L^{\theta}\circ f_{L-1}^{\theta}\circ...\circ f_1^{\theta})(z),
    \end{equation*}
    where in the setting of the semilinear wave equation on $\Omega\subset\mathbb{R}^d, d\in\mathbb{N}$ we set $l_0=d+1$ and $z=(x,t)$. We refer to $u_{\theta}$ as the realization of the \textit{neural network} associated with the parameters $\theta$ with $L$ layers and widths $(l_0,l_1,...,l_L)$. We refer to layers $1$ to $L-1$ as \textit{hidden layers}. For $1\leq K\leq L$, we say that layer $k$ has width $l_k$, and we refer to $W_k$ and $b_k$ as the \textit{weights} and \textit{biases} corresponding to layer $k$. The width $W$ of $u_{\theta}$ is defined as $\mathrm{max}(l_0,...,l_L)$. If $L=2$, we that that $u_{\theta}$ is a \textit{shallow neural network}; if $L\geq 3$, then we say $u_{\theta}$ is a \textit{deep neural network}. 
\end{definition}

Figure \ref{fig:NN} illustrates the architecture of neural network with two layers employed for the semilinear wave equation.  

\begin{figure}
    \centering

\begin{tikzpicture}[>=stealth, node distance=1.5cm]

\foreach \name / \y / \text in {1/1/x,2/0/y,3/-1/t}
    \node[draw, circle, draw=green, fill=green!20, minimum size=22] (I-\name) at (0, \y) {\text};
\node[above=0.5cm of I-1, green] {Input Layer};

\foreach \name / \y in {1/-2,2/-1,3/0,4/1,5/2}{
    \ifnum \y=0
         \node at (2,\y) {$\vdots$};
    \else
        \node[draw, circle, draw = blue, fill=blue!20, minimum size=22] (H1-\name) at (2, \y) {};
    \fi}
\node[above =0.4cm of H1-5, xshift=1cm, blue] {Hidden Layers};

\foreach \name / \y in {1/-2,2/-1,3/0,4/1,5/2}{
    \ifnum \y=0
         \node at (4,\y) {$\vdots$};
    \else
         \node[draw, circle, draw = blue, fill=blue!20, minimum size=22] (H2-\name) at (4, \y) {};
    \fi}

\node[draw, circle, draw = red, fill=red!20, minimum size=22] (O-1) at (6, 0) {$u_{\theta}(x, y, t)$};
\node[above=0.5cm of O-1, red] {Output Layer};

\foreach \source in {1,2,3}
    \foreach \dest in {1,2,4,5}
        \draw[->] (I-\source) -- (H1-\dest);

\foreach \source in {1,2,4,5}
    \foreach \dest in {1,2,4,5}
        \draw[->] (H1-\source) -- (H2-\dest);

\foreach \source in {1,2,4,5}
    \draw[->] (H2-\source) -- (O-1);

\end{tikzpicture}
\caption{Architecture of the neural network with three input nodes, two hidden layers with 80 neurons each and one output node.}
\label{fig:NN}
\end{figure}

Moreover, the midpoint rule will also be needed throughout:
\begin{definition}[Midpoint rule \cite{DeRyck1}]
    Given $\Omega\subset\mathbb{R}^d$ and $f\in L^1(\Omega)$ partition $\Omega$ into $M\sim N^d$ cubes of edge length $\frac{1}{N}$ for $N\in\mathbb{N}$ and denote by $\{y_m\}_{m=1}^M$ the midpoint of these cubes. Then the midpoint rule $\mathcal{Q}_M^{\Omega}$ and its accuracy are given by
\begin{equation}\label{Quadrature2}
    \begin{split}
        &\mathcal{Q}_M^{\Omega}[f]:=\frac{1}{M}\sum_{m=1}^Mf(y_m),\\
        &|\int_{\Omega}f(y)dy-\mathcal{Q}_M^{\Omega}[f]| \leq C_{\Omega}\|f\|_{C^2}M^{-\frac{2}{d}},
    \end{split}
\end{equation}
where $C_{\Omega}$ depends on the the geometry of $\Omega$.
\end{definition}
\section{The Semilinear Wave Equation}
The following sections are devoted to the main result of the work. To that end, first the semilinear wave equation will be introduced, as well as corresponding existence results for an analytic solution. Then the three questions proposed in \cite{DeRyck1} will be answered consecutively. In the end, the main result consists of an \textit{a posteriori} and \textit{a priori} error bound for the sum of the total error and the total error of the time-derivative in terms of the training error, the number of training points and some constants. Some of the longer proofs of the theorems and corollaries can be found in the appendix.\\
The focus of this work will be the following semilinear wave equation:
\begin{equation}\label{WaveBVP}
\begin{cases}
        u_{tt}-\Delta u+a(x)u_t+f(x,u)=0\qquad &\text{on }\Omega\times[0,T]\\
        u(x,0)=u_0(x)&\text{on }\Omega\\
        u_t(x,0)=u_1(x)&\text{on }\Omega\\
        u(x,t)=0&\text{on }\partial\Omega\times[0,T],
\end{cases}
\end{equation}
where $\Omega$ is a bounded domain in $\mathbb{R}^d$ and $d\in\mathbb{N}$ with $d\leq 5$ and $a:\overline{\Omega}\rightarrow\mathbb{R}, f:\overline{\Omega}\times\mathbb{R}\rightarrow\mathbb{R}, u_0:\Omega\rightarrow\mathbb{R}$ and $u_1:\Omega\rightarrow\mathbb{R}$. In practice, the nonlinearity $f(x,u)$ often takes the shape $f(x,u)=|u|^pu$ for some $p>0$.\\
Throughout we will make the following assumption on $a(x)$ and $f(x,u)$:
\begin{itemize}\label{Assumptions}
    \item[(A1)] The function $a:\overline{\Omega}\rightarrow\mathbb{R}$ satisfies $a(x)\geq 0, a(x)\neq 0$ almost everywhere in $\overline{\Omega}$ and $a(x)^{-1}\equiv\frac{1}{a(x)}\in L^p(\Omega)$ for some $0<p<\infty$.
    \item[(A2)]\label{A2} The function $f(x,u)$ is a $C^{k-1}$ class function defined on $\overline{\Omega}\times\mathbb{R}$ satisfying
    \begin{equation}
        \left| D^i \left( \frac{\partial}{\partial u} \right) ^jf(x,u) \right|\leq c|u|^{(r+1-j)^+}\qquad 0\leq i+j\leq k-1,\forall x\in\overline{\Omega}, \forall u\in\mathbb{R},
    \end{equation}
    where $c>0, k\in\mathbb{N},k>3$ and $D^i$ denotes the partial differentiations of order $i$ with respect to $x=(x_1,...,x_d)$ and $r$ satisfies
    \begin{equation}
    \begin{split}
        &1\leq r<\infty\qquad\qquad\text{if }d\leq 2,\\
        &1\leq r< \frac{d+2}{d-2}\qquad\text{if }2< d\leq 4,\\
        &1\leq r<2\qquad\qquad\text{if }d=5.
    \end{split}
    \end{equation}
\end{itemize}
Our error estimation results rely on the following existence, uniqueness and regularity results that can be found in \cite{Nakao}. The first result guarantees local existence for some fixed $T>0$, whereas the second result yields global existence for all $T$.
\begin{proposition}[Local Existence]\label{LocalExistence}
    Let $1\leq d<6$ and assume $\Omega\subset\mathbb{R}^d$ has a $C^3$ boundary. Suppose that (A1) and (A2) are fulfilled with $k\in\mathbb{N}, k>3$. Let $(u_0,u_1)\in D((-\Delta)^{\frac{k+1}{2}})\times D((-\Delta)^{\frac{k}{2}})$ where $D(-\Delta)=H^2(\Omega)\cap H^1_0(\Omega)$. Then there exists $T=T(\|u_0\|_{H^{k+1}(\Omega)},\|u_1\|_{H^k(\Omega)})>0$ such that the problem \ref{WaveBVP} admits a unique solution $u$ in the class
        \begin{equation*}
        C^{k+1}([0,T);L^2(\Omega))\bigcap_{i=0}^kC^i([0,T);H^{k+1-i}(\Omega)\cap H^1_0(\Omega))
    \end{equation*}
\end{proposition}
\begin{theorem}[Global Existence]\label{GlobalExistence}
    Let $1\leq d<6$ and assume $\Omega\subset\mathbb{R}^d$ has a $C^3$ boundary. Suppose that (A1) and (A2) are fulfilled with $k\in\mathbb{N}, k>3$ and $p>p_0(d,r)$ with
    \begin{equation*}
        p_0=
        \begin{cases}
            \frac{d}{3r} & \text{if }d\leq 2,\\
            \frac{2d}{(8-d)r+d-2} & \text{if }2<s\leq 4\text{ and }1\leq r<\frac{2}{d-2},\\
            \mathrm{max}\{\frac{2d}{(8-d)r+d-2},\frac{2d}{3(4-d)r+6}\}& \text{if }2<s\leq 4\text{ and }\frac{2}{d-2}\leq r<\frac{2d-2}{d-2},\\
            \mathrm{max}\{\frac{2d}{3(4-d)r+5d-6},\frac{2d}{3(4-d)r+6}\}& \text{if }2<s\leq 4\text{ and }\frac{2d-2}{d-2}\leq r,\\
            \frac{10}{6-3r}&\text{if }d=5
        \end{cases}
    \end{equation*}
    Then there exists an open unbounded set $\mathcal{U}\subset D((-\Delta)^2)\times D((-\Delta)^{\frac{3}{2}})$ such that if $(u_0,u_1)\in\mathcal{U}\cap D((-\Delta)^{\frac{k+1}{2}})\times D((-\Delta)^{\frac{k}{2}})$ the problem \ref{WaveBVP} admits a unique global $(T=\infty)$ solution $u$ in the class
    \begin{equation*}
        C^{k+1}(\mathbb{R}^+;L^2(\Omega))\bigcap_{i=0}^kC^i(\mathbb{R}^+;H^{k+1-i}(\Omega)\cap H^1_0(\Omega)).
    \end{equation*}
    Moreover, the solutions is a classical solution.
\end{theorem}
We note that the set $\mathcal{U}$ contains a neighbourhood of the origin in $D((-\Delta)^2)\times D((-\Delta)^{\frac{3}{2}})$ and depends in a non-trivial way on the norms of the local solution that we get from \ref{LocalExistence} (details in \cite{Nakao}). Moreover, the regularity of the solution depends on the regularity of the initial data. Therefore, higher regularity implies faster convergence rates, see Theorem \ref{Q1} Corollaries \ref{BoundGeneralizationError} and \ref{BountTraininError}.

Based on this result, we can prove that $u$ is Sobolov-regular, i.e. $u\in H^{k+1}(\Omega\times[0,T])$ for some $k\in\mathbb{N}$. 
\begin{corollary}\label{Regularity}
    Under the assumptions of \ref{GlobalExistence}, for every $T>0$ there exists a unique solution $u$ to the semi-linear wave equation such that $u\in H^{k+1}(\Omega\times(0,T))$.
\end{corollary}
\begin{proof}
    By theorem \ref{GlobalExistence}, we know that for every $T>0$ there exists a unique classical solution $u$ to the semilinear wave equation \ref{WaveBVP} in the class 
    \begin{equation*}
        C^{k+1}([0,T];L^2(\Omega))\bigcap_{i=0}^kC^i([0,T];H_{k+1-i}(\Omega)\cap H^1_0(\Omega)).
    \end{equation*}
    To show that $u\in H^{k+1}(\Omega\times(0,T))$ we need to show that all weak derivatives of $u$ up to order $k+1$ exist and are in $L^2(\Omega\times(0,T))$. But this follows immediately from the fact that
    \begin{equation*}
        C([0,T];X)\subset L^p(0,T;Y)
    \end{equation*} for all $1\leq p\leq\infty$ and Banach spaces $X,Y$ with $X\subset Y$. 
\end{proof}
\section{Error Estimation for Semilinear Wave Equation}
\subsection{PINN Residuals}
To investigate the PINN solution to \ref{WaveBVP}, we consider the pointwise residuals, defined for the semilinear wave equation \ref{WaveBVP} and for any sufficiently smooth function $v:\Omega\times[0,T]\rightarrow\mathbb{R}$ as
\begin{equation}\label{Residuals1}
    \begin{split}
    &\mathcal{R}_{PDE}[v](x,t)=v_{tt}-\Delta v+a(x)v_t+f(x,v),\,\qquad\,\,(x,t)\in\Omega\times[0,T],\\
    &\mathcal{R}_{s,u}[v](y,t)=v(y,t),\quad\qquad\qquad\qquad\qquad\qquad\quad\, (x,t)\in\partial\Omega\times[0,T],\\
    &\mathcal{R}_{u_0}[v](x)=v(x,0)-u_0,\qquad\qquad\qquad\qquad\qquad\quad \forall x\in\Omega,\\
    &\mathcal{R}_{u_1}[v](x)=v_t(x,0)-u_1,\qquad\qquad\qquad\qquad\qquad\:\,\ \forall x\in\Omega.\\
    \end{split}
\end{equation}
The above residuals can be considered the canonical residuals that enforce the boundary conditions on the boundary $\partial\Omega\times[0,T]$ ($\mathcal{R}_{s,u}[v])$, the initial conditions on $\Omega$ ($\mathcal{R}_{u_0}[v], \mathcal{R}_{u_1}[v]$) and the PDE dynamics in $\Omega\times[0,T]$ ($\mathcal{R}_{PDE}[v]$). In addition, for practical reasons that will become clear later when proving Theorem \ref{Q1}, we also need the following two additional residuals:
\begin{equation}\label{Residuals2}
    \begin{split}
        &\mathcal{R}_{s,u_t}[v](y,t)=v_t(y,t)\,\,\qquad\qquad\:(x,t)\in\partial\Omega\times[0,T],\\
        &\mathcal{R}_{\nabla u}[v](x)=\nabla v(x,0)-\nabla u_0\qquad\forall x\in\Omega,
    \end{split}
\end{equation}
where $\mathcal{R}_{s,u_t}[v]$ is motivated by the fact that from $u(x,t)=0$ for all $(y,t)\in\partial\Omega\times[0,T]$ it follows that $u_t(y,t)=0$ for all $(y,t)\in\partial\Omega\times[0,T]$. For the exact solution $u$ of \ref{WaveBVP} it holds that $\mathcal{R}_{PDE}[u]=\mathcal{R}_{s,u}[u]=\mathcal{R}_{u_0}[u](x)=\mathcal{R}_{u_1}[u]=\mathcal{R}_{s,u_t}[u]=\mathcal{R}_{\nabla u}[u]=0$.

The PINN algorithm aims to find a neural network $u_{\theta}$ that minimizes all residuals simultaneously. This is done by minimizing the \textit{generalization error}
\begin{equation}\label{GenError}
\begin{split}
        \mathcal{E}_G^2(\theta)=&\int_{\Omega\times[0,T]}|\mathcal{R}_{PDE}[u_{\theta}](x,t)|^2dxdt+\int_{\partial\Omega\times[0,T]}|\mathcal{R}_{s,u}[u_{\theta}](s,t)|^2ds(x)dt\\
        &+\int_{\partial\Omega\times[0,T]}|\mathcal{R}_{s,u_t}[u_{\theta}](s,t)|^2ds(x)dt+\int_{\Omega}|\mathcal{R}_{u_0}[u_{\theta}](x)|^2dx+\int_{\Omega}|\mathcal{R}_{u_1}[u_{\theta}](x)|^2dx\\
        &+\int_{\Omega}\|\mathcal{R}_{\nabla u}[u_{\theta}](x)\|^2_{\mathbb{R}^d}dx.
\end{split}
\end{equation}
The generalization error involves the minimization of integrals which cannot be directly computed in practice, which is why a quadrature rule is used instead to approximate these integrals. This gives rise to the \textit{training error}
\begin{equation}\label{TrainingErr}
    \begin{split}
        \mathcal{E}_T(\theta,\mathcal{S})^2&=\mathcal{E}_T^{PDE}(\theta,\mathcal{S}_{PDE})^2+\mathcal{E}_T^{s,u}(\theta,\mathcal{S}_{s})^2+\mathcal{E}_T^{s,u_t}(\theta,\mathcal{S}_{s})^2+\mathcal{E}_T^{u_0}(\theta,\mathcal{S}_{t})^2+\mathcal{E}_T^{u_1}(\theta,\mathcal{S}_{t})^2+\mathcal{E}_T^{\nabla u}(\theta,\mathcal{S}_{t})^2\\
        &=\sum_{n=1}^{M_{PDE}}w_{PDE}^n|\mathcal{R}_{PDE}[u_{\theta}](x_{PDE}^n,t_{PDE}^n)|^2+\sum_{n=1}^{M_s}w_{s,u}^n|\mathcal{R}_{s,u}[u_{\theta}](x_s^n,t_s^n)|^2\\
        &\quad+\sum_{n=1}^{M_s}w_{s,u_t}^n|\mathcal{R}_{s,u_t}[u_{\theta}](x_s^n,t_s^n)|^2+\sum_{n=1}^{M_t}w_{u_0}^n|\mathcal{R}_{u_0}[u_{\theta}](x_t^n)|^2\\
        &\quad+\sum_{n=1}^{M_t}w_{u_1}^n|\mathcal{R}_{u_1}[u_{\theta}](x_t^n)|^2+\sum_{n=1}^{M_t}w_{\nabla u}^n\|\mathcal{R}_{\nabla u}[u_{\theta}](x_t^n)\|^2_{\mathbb{R}^d}
    \end{split}
\end{equation}
with quadrature points 
\begin{equation*}
    \mathcal{S}_{PDE}=\{(x_{PDE}^n,t_{PDE}^n)\}_{n=1}^{M_{PDE}},\,\mathcal{S}_s=\{(x_s^n,t_s^n)\}_{n=1}^{M_s},\,\mathcal{S}_t=\{(x_t^n)\}_{n=1}^{M_t},
\end{equation*}
$M_{PDE}, M_s, M_t\in\mathbb{N}$ and the corresponding quadrature weights $w_{PDE}^n, w_{s,u}^n,w_{s,u_t}^n, w_{u_0}^n, w_{u_1}^n, w_{\nabla u}^n\in\mathbb{R}^+$.

The PINN algorithm finds an approximation $u_{\theta}$ to \ref{WaveBVP} by attempting to solve the following optimization problem
\begin{equation}\label{OptimizationProblem2}
    \theta^*=\underset{\theta\in\Theta}{arg\,min\,}\mathcal{E}_T(\theta,\mathcal{S})^2.
\end{equation}
In general, we have to keep in mind that there is no guarantee that the algorithm finds the global minimum, since the problem is often highly non-convex.
\subsection{Question 1}
Question 1 asks whether there exists a neural network $u_{\theta}$, such that the associated generalization error and training error are arbitrarily small. The following theorem gives an affirmative answer to this question and is based on a result in \cite{DeRyck1} that bounds the $H^k$-norms of the difference between a function $u$ and its neural network approximation $u_{\theta}$ for neural networks with a specific architecture. Theorem \ref{Q1} has two parts: The first is the existence and uniqueness result for the semilinear wave equation and the second part are the bounds for the individual residuals. The existence result is needed to guarantee that the algorithm has something it can converge to in the first place, i.e. an existing actual solution of the semi-linear wave equation. 
\begin{theorem}\label{Q1}
    Suppose the assumptions of \ref{GlobalExistence} are fulfilled (notably $k\in\mathbb{N},k> 3)$ and assume additionally that $a(x)\in L^{\infty}(\Omega)$ and that there exists $\gamma\in[0,1)$ such that
    \begin{equation}\label{r}
        r=\frac{4\gamma}{d+1-4\gamma}.
    \end{equation}
    Let $N\in\mathbb{N}, n\geq2, \delta>0$ and let $\Omega\subset\mathbb{R}^d, 1\leq d<6$, have a $C^3$ boundary. Let $\hat{\Omega}$ denote a subset of $\mathbb{R}^d$ such that $\Omega\subset\subset\hat{\Omega}$ and $\hat{\Omega}=\Pi_{i=1}^d[a_i,b_i]$ for some $a_i,b_i\in\mathbb{Z},a_i<b_i$. We use the notation $\Omega_T:=\Omega\times[0,T]$ and $\hat{\Omega}_T:=\hat{\Omega}\times[0,T]$. Then it holds that 
    \begin{enumerate}[I.]
        \item for all $T>0$ there exists a unique classical solution $u$ to the semi-linear wave equation \ref{WaveBVP} such that $u\in H^{k+1}(\Omega\times(0,T))$.
        \item For every $N>5$, there exists a tanh neural network $u_{\theta}$ with two hidden layers of width $3\lceil \frac{k+n-1}{2}\rceil\binom{d+k+1}{k}+\lceil (N-1)(T+\sum_{i=1}^d(b_i-a_i) \rceil$ and $3(d+3)N^{d+1}\lceil \frac{d+n+1}{2} \rceil \lceil T\Pi_{i=1}^d(b_i-a_i) \rceil$ such that 
        \begin{equation*}
        \begin{split}
             &\|\mathcal{R}_{PDE}[u_{\theta}]\|_{L^2(\Omega_T)}=\|u_{\theta_{tt}}-\Delta u_{\theta}+a(x)u_{\theta_t}+f(x,u_{\theta})\|_{L^2(\Omega_T)}\\
             &\leq (C_2+\sqrt{d})\lambda_2(N)N^{-k+1}+\|a\|_{L^{\infty}(\Omega)}C_1\lambda_1(N)N^{-k}+\tilde{C}(C_2\lambda_2(N)N^{-k+1})^{\gamma}(C_0\lambda_0(N)N^{-k-1})^{1-\gamma},\\
             &\|\mathcal{R}_{s,u}[u_{\theta}]\|_{L^2(\partial\Omega\times[0,T])}=
             \|u_{\theta}\|_{L^2(\partial\Omega\times[0,T])}\leq C_{\Omega_T}C_1\lambda_1(N)N^{-k},\\
             &\|\mathcal{R}_{s,u_t}[u_{\theta}]\|_{L^2(\partial\Omega\times[0,T])}=
             \|u_{\theta_t}\|_{L^2(\partial\Omega\times[0,T])}\leq C_{\Omega_T}C_2\lambda_2(N)N^{-k+1},\\
             &\|\mathcal{R}_{u_0}[u_{\theta}]\|_{L^2(\Omega)}=
             \|u_{\theta}(t=0)-u_0\|_{L^2(\Omega)}\leq C_{\Omega_T}C_1\lambda_1(N)N^{-k},\\
             &\|\mathcal{R}_{u_1}[u_{\theta}]\|_{L^2(\Omega)}=
             \|u_{\theta_t}(t=0)-u_1\|_{L^2(\Omega)}\leq C_{\Omega_T}C_2\lambda_2(N)N^{-k+1},\\
             &\|\mathcal{R}_{\nabla u}[u_{\theta}]\|_{L^2(\Omega)}=
             \|\nabla u_{\theta}(t=0)-\nabla u_0\|_{L^2(\Omega)}\leq \sqrt{d}C_{\Omega_T}C_2\lambda_2(N)N^{-k+1},
        \end{split}
        \end{equation*}
        where $\tilde{C}$ is defined in the proof and depends on $N, \|u\|_{H^2(\Omega_T)}$ and $\|u\|_{L^2(\Omega_T)}$ and
        \begin{equation*}
            \begin{split}
                &\lambda_l(N)=2^l3^{d+1}(1+\delta)\mathrm{ln}^l(\beta_{l,\delta}N^{d+k+4}),\\
                &\beta_{l,\delta}=\frac{2^{l(d+1)}5\mathrm{max}\{T\Pi_{i=1}^d(b_i-a_i),d+1\}\mathrm{max}\{\|u-u_{\theta}\|_{W^{l,\infty}(\Omega_T)},1\}}{3^{d+1}\delta\mathrm{min}\{1,C_l\}},\\
                &C_l=\underset{0\leq i\leq l}{\mathrm{max}}\binom{d+i}{i}^{\frac{1}{2}}\frac{((k+1-i)!)^{\frac{1}{2}}}{(\lceil\frac{k+1-i}{d+1}\rceil !)^{\frac{d+1}{2}}}\left(\frac{3\sqrt{d+1}}{\pi}\right)^{k+1-i}|u-u_{\theta}|_{H^{k+1}(\Omega_T)},\\
                &C_{\Omega_T}=\sqrt{\frac{2\mathrm{max}\{2h_{\Omega_T},d+1\}}{\rho_{\Omega_T}}},
            \end{split}
        \end{equation*}
        where $h_{\Omega_T}$ denotes the diameter of $\Omega_T$ and $\rho_{\Omega_T}$ the radius of the largest $d+1$-dimensional ball that can be inscribed into $\Omega_T$. The weights of the network can be bounded by $\mathcal{O}(N\mathrm{ln}(N)+N^{\kappa})$, where $\kappa=\mathrm{max}\{(k+1)^2,(d+1)(d+k+4)\}/n$.
    \end{enumerate}
\end{theorem}

This result shows that we can derive upper bounds for all the individual residuals and that these decrease with $N$. Moreover, the weights of the network grow with $\mathcal{O}(N\mathrm{ln}(N)+N^{\kappa})$, where $\kappa=\mathrm{max}\{(k+1)^2,(d+1)(d+k+4)\}/n$. Therefore, if $n$ is big enough such that $\kappa\leq 1$ the weights only grow with $\mathcal{O}(N\mathrm{ln}(N))$ which becomes important for the question whether the convergence of the training error implies the convergence of the total error (see Theorem \ref{apriori}). Both observations would imply that bigger $N$ and $n$ are beneficial for getting smaller errors. However, this comes at a cost in the form of the width of the network which grows in both $n$ and $N$. Nevertheless, it is an interesting aspect of this result that to get better error bounds we do not need more than two hidden layers, instead focusing on the width of the network. This is different from the usual paradigm in machine learning where better performance is often associated with more depth.

The upper bounds on the individual residuals that we have computed in Theorem \ref{Q1} directly give us an upper bound on the generalization error. 
\begin{corollary}\label{BoundGeneralizationError}
    Under the same assumptions as Theorem \ref{Q1}, it follows directly from Theorem \ref{Q1} that there exists a neural network $u_{\theta}$ whose generalization error can be bounded in the following way and can be made arbitrarily small. Let $S:=\{PDE, (s,u), (s,u_t), u_0, u_1, \nabla u\}$, then
    \begin{equation*}
        \begin{split}
            \mathcal{E}_G(\theta)&=\sqrt{\sum_{i\in S}\|\mathcal{R}_i\|^2}\leq\sum_{i\in S}\|\mathcal{R}_i\|=(C_2(1+C_{\Omega_T}(2+\sqrt{d}))+\sqrt{d})\lambda_2(N)N^{-k+1}\\
            &\quad+(\|a\|_{L^{\infty}(\Omega)}+2C_{\Omega_T})\lambda_1(N)N^{-k}+\tilde{C}C_2^{\gamma}C_0^{1-\gamma}\lambda_2(N)^{\gamma}\lambda_0(N)^{1-\gamma}N^{-k-1+2\gamma}\\
            &=(C_2(1+C_{\Omega_T}(2+\sqrt{d}))+\sqrt{d})2^23^{d+1}(1+\delta)\mathrm{ln}^2(\beta_{2,\delta}N^{d+k+4})N^{-k+1}\\
            &\quad+(\|a\|_{L^{\infty}(\Omega)}+2C_{\Omega_T})3^{d+1}2(1+\delta)\mathrm{ln}(\beta_{1,\delta}N^{d+k+4})N^{-k}\\
            &\quad+\tilde{C}C_2^{\gamma}C_0^{1-\gamma}(2^23^{d+1}(1+\delta)\mathrm{ln}^2(\beta_{2,\delta}N^{d+k+4}))^{\gamma}(3^{d+1}(1+\delta))^{1-\gamma}N^{-k-1+2\gamma}\\
            &=\mathcal{O}(\mathrm{ln^2(N)N^{-k+1}}),
        \end{split}
    \end{equation*}
    which goes to $0$ for $N\rightarrow\infty$.
\end{corollary}
Next, we use Theorem \ref{Q1} and the quadrature rule \ref{Quadrature2} to find an upper bound for the training error of the neural network we have constructed in Theorem \ref{Q1}, thus showing that there exists a physics-informed neural network with arbitrarily small training error. For simplicity's sake, the training set is obtained with the midpoint rule $\mathcal{Q}_M$ \ref{Quadrature2}. Other quadrature rules could be used as well to achieve a slightly better convergence rate but might incur higher computational costs. For higher dimensions, numerical quadrature rules are usually considered not suitable and instead random sampling techniques, such as Monte-Carlo integration, are used instead \cite{Mishra1}. To apply the quadrature rule \ref{Quadrature2} we assume that for an $\Lambda\subset\mathbb{R}^d$ there exists an $N\in\mathbb{N}$ and a partition of $\Lambda$ into $M\sim N^d$ cubes of side-length $\frac{1}{N}$ such that the midpoints of these cubes are inside $\Lambda$. Then for $f=\mathcal{R}_{PDE}^2$ and $\Lambda=\Omega\times(0,T)$ we obtain the quadrature $\mathcal{Q}_{M_{PDE}}^{PDE}$, for $f\in\{\mathcal{R}_{u_0}^2,\mathcal{R}_{u_1}^2, \mathcal{R}_{\nabla u_0}^2\}$ and $\Lambda=\partial\Omega\times[0,T]$ we obtain the quadrature $\mathcal{Q}_{M_t}^t$ and for $f\in\{\mathcal{R}_{s,u}^2,\mathcal{R}_{s,u_t}^2\}$ and $\Lambda = \Omega$ we obtain the quadrature $\mathcal{Q}_{M_s}^s$.

The training error \ref{TrainingErr} can then be written as follows:
\begin{equation}\label{ErrorQuad}
    \begin{split}
        \mathcal{E}_T(\theta,\mathcal{S})^2&=\mathcal{E}_T^{PDE}(\theta,\mathcal{S}_{PDE})^2+\mathcal{E}_T^{s,u}(\theta,\mathcal{S}_{s})^2+\mathcal{E}_T^{s,u_t}(\theta,\mathcal{S}_{s})^2+\mathcal{E}_T^{u_0}(\theta,\mathcal{S}_{t})^2+\mathcal{E}_T^{u_1}(\theta,\mathcal{S}_{t})^2+\mathcal{E}_T^{\nabla u}(\theta,\mathcal{S}_{t})^2\\
        &=\mathcal{Q}^{PDE}_{M_{PDE}}[\mathcal{R}_{PDE}^2]+\mathcal{Q}^{s,u}_{M_s}[\mathcal{R}_{s,u}^2]+\mathcal{Q}^{s,u_t}_{M_s}[\mathcal{R}_{s,u_t}^2]+\mathcal{Q}^{u_0}_{M_t}[\mathcal{R}_{u_0}^2]+\mathcal{Q}^{u_1}_{M_t}[\mathcal{R}_{u_1}^2]+\mathcal{Q}^{\nabla u}_{M_t}[\mathcal{R}_{\nabla u}^2],
    \end{split}
\end{equation}
where the weights in \ref{TrainingErr} are the quadrature weights of the midpoint rule $w_{PDE}^n=\frac{1}{M_{PDE}}, w_{s,u}^n,w_{s,u_t}^n=\frac{1}{M_s}$ and $w_{u_0}^n, w_{u_1}^n, w_{\nabla u}^n=\frac{1}{M_t}$ and the training sets $\mathcal{S}_{PDE}=\{(x_{PDE}^n,t_{PDE}^n)\}_{n=1}^{M_{PDE}}, \mathcal{S}_s=\{(x_s^n,t_s^n)\}_{n=1}^{M_s}, \mathcal{S}_t=\{(x_t^n)\}_{n=1}^{M_t}$ consist of the midpoints of the cubes.
\begin{corollary}\label{BountTraininError}
    Let the assumptions of Theorem \ref{Q1} be fulfilled and let let $\Omega\subset\mathbb{R}^d$ be such that there exists an $\tilde{N}\in\mathbb{N}$ and a partition of $\Omega$ into $M\sim \tilde{N}^d$ cubes of side-length $\frac{1}{\tilde{N}}$ such that the midpoints of these cubes are inside $\Omega$. Then for every $N\in\mathbb{N}, N>5$ there exists a tanh neural network $u_{\theta}$ with two hidden layers of width $3\lceil \frac{k+n-1}{2}\rceil\binom{d+k+1}{k}+\lceil (N-1)(T+\sum_{i=1}^d(b_i-a_i) \rceil$ and $3(d+3)N^{d+1}\lceil \frac{d+n+1}{2} \rceil \lceil T\Pi_{i=1}^d(b_i-a_i) \rceil$ such that we get the following result for the training error:
    \begin{equation*}
        \mathcal{E}_T(\theta,\mathcal{S})^2=\mathcal{O}(\mathrm{ln}^4(N)N^{2(-k+1)}+M^{-\frac{2}{d}}_s-M^{-\frac{2}{d}}_t+M^{-\frac{2}{d+1}}_{PDE})
    \end{equation*}
\end{corollary}
\begin{proof}
    Theorem \ref{Q1} gives us the existence of a classical solution to the semi-linear wave equation \ref{WaveBVP}, as well as the existence of a tanh neural network such that all residuals are bounded and according to Corollary \ref{BoundGeneralizationError} the generalization error grows like $\mathcal{O}(\mathrm{ln^2(N)N^{-k+1}})$. From the quadrature rule we know that
    \begin{equation*}
        \mathcal{Q}^{\Omega}_M[f]\leq C_{\Omega}\|f\|_{C^2}M^{-\frac{2}{d}}+\int_{\Omega}fd(y).
    \end{equation*}
    These two results lead directly to the following bound for the training error:
    \begin{equation*}
        \begin{split}
            \mathcal{E}_T(\theta,\mathcal{S})^2&\leq 4C_{\Omega_T}\|\mathcal{R}_{PDE}\|_{C^2}^2M^{-\frac{2}{d+1}}_{PDE}+\|\mathcal{R}_{PDE}\|_{L^2(\Omega_T)}^2+ 4C_{\partial\Omega\times[0,T]}\|\mathcal{R}_{s,u}\|_{C^2}^2M^{-\frac{2}{d}}_s+\|\mathcal{R}_{s,u}\|_{L^2(\partial\Omega\times[0,T])}^2\\
            &\quad+4C_{\partial\Omega\times[0,T]}\|\mathcal{R}_{s,u_t}\|_{C^2}^2M^{-\frac{2}{d}}_s+\|\mathcal{R}_{s,u_t}\|_{L^2(\partial\Omega\times[0,T])}^2+4C_{\Omega}\|\mathcal{R}_{u_0}\|_{C^2}^2M^{-\frac{2}{d}}_t+\|\mathcal{R}_{u_0}\|_{L^2(\Omega)}^2\\
            &\quad+4C_{\Omega}\|\mathcal{R}_{u_1}\|_{C^2}^2M^{-\frac{2}{d}}_t+\|\mathcal{R}_{u_1}\|_{L^2(\Omega)}^2+4C_{\Omega}\|\mathcal{R}_{\nabla u}\|_{C^2}^2M^{-\frac{2}{d}}_t+\|\mathcal{R}_{\nabla z}\|_{L^2(\Omega)}^2\\
            &=\mathcal{O}(\mathcal{E}_G(\theta)^2+M^{-\frac{2}{d}}_s+M^{-\frac{2}{d}}_t+M^{-\frac{2}{d+1}}_{PDE})\\
            &=\mathcal{O}(\mathrm{ln}^4(N)N^{2(-k+1)}+M^{-\frac{2}{d}}_s-M^{-\frac{2}{d}}_t+M^{-\frac{2}{d+1}}_{PDE}),
        \end{split}
    \end{equation*}
    where we have used that for every residual it holds that $\|\mathcal{R}_i^2\|_{C^n}\leq 2^n\|\mathcal{R}_i\|^2_{C^n}$ (from the general Leibnitz rule).
\end{proof}
This result shows that there exists a neural network $u_{\theta}$ such that for $N,M$ large enough, the training error can be made arbitrarily small. However, this does however not guarantee that we can find such a PINN in practice.
\begin{remark}
The assumptions that $\Omega\subset\mathbb{R}^d$ has a $C^3$ boundary and that $1\leq d<6$ in Theorem \ref{Q1} are necessary to guarantee the existence of a classical solution to the semilinear wave equation \ref{WaveBVP} in accordance with Corollary \ref{Regularity}. This is needed to make sure that there exists a function that the algorithm can approximate in the first place. If the existence of a sufficiently regular solution is known by other means, for example by explicit computation or a different existence theorem, then these requirements can be dropped while the result stays the same. For this reason, these requirements will be dropped for question 2 and parts of question 3.
\end{remark}
\begin{remark}
A tanh neural network was chosen to ensure the applicability of Theorem \ref{tanh}. However, remarks in \cite{DeRyck2} imply that this result also holds for other smooth activation functions such as logistic or sigmoid.
\end{remark}
\subsection{Question 2}
Question 2 considers whether the total error for a PINN $u_{\theta}$ is small if the generalization error is also small. To this end we bound the $L^2$-norm of the total error by the $L^2$-norm of the residuals times a constant. Moreover, it turns out that we can not only bound the total error but also the total error of the time derivative. Both can be made small if the generalization error, i.e. the residuals, are small. Importantly, this result is independent of the particular neural network architecture that we have used to derive an affirmative answer to question 1.
\begin{theorem}\label{Q2}
    Let $T>0, d\in\mathbb{N}, \Omega\subset\mathbb{R}^d$ be a bounded Lipschitz domain and let $u\in C^2(\Omega\times[0,T])$ be a classical solution of the semi-linear wave equation \ref{WaveBVP} and let the assumptions (A1) and (A2) hold. Let $u_{\theta}$ be a PINN with parameters $\theta$ and let $\Omega_T:=\Omega\times[0,T]$. Then we can bound the total error and the error of the time-derivative in the $L^2$-norm:
    \begin{equation*}
        \int_{\Omega_T}|u_{\theta}-u|^2dxdt+\int_{\Omega_T}|u_{\theta_t}-u_t|^2dxdt\leq\mathcal{C}T\mathrm{exp}\left(T(\mathrm{max}\{1,2C_{pw}^2\}(1+\hat{C}+\frac{2\sqrt{T}}{C_{pw}^2}))\right),
    \end{equation*}
    where $\hat{C}$ is defined in the proof and depends on $\|u\|_{C^0},\|u_{\theta}-u\|_{C^0}$ and $r$, $C_{pw}$ is the constant from an application of the Poincar\'e-Wirtinger inequality, see \cite{evans}, and
    \begin{equation*}
    \begin{split}
        \mathcal{C}&=\text{max}\{1,2C^2_{pw}\}\biggl(\|\mathcal{R}_{u_1}\|_{L^2(\Omega)}^2+\|\mathcal{R}_{PDE}\|_{L^2(\Omega_T)}^2+\|\mathcal{R}_{\nabla u}\|^2_{L^2(\Omega)}\\
        &\quad+2\sqrt{T|\partial\Omega|}\|u_\theta-u\|_{C^1}\|\mathcal{R}_{s,u_t}\|_{L^2(\partial\Omega\times[0,T])}+\frac{2}{C_{pw}^2}\|\mathcal{R}_{u_0}\|_{L^2(\Omega)}^2\biggr).
    \end{split}
    \end{equation*}
\end{theorem}

This result shows that the total error of a physics-informed neural network can be bound in terms of the residuals. For neural networks that have a small generalization error, such as the one from Theorem \ref{Q1}, it is therefore guaranteed that the total error will be small as well. Moreover, Theorem \ref{Q2} goes above just showing a bound on the total error but gives a stronger result: both the total error and the error of the time-derivative are bounded and can be made small if the generalization error is small. This suggests that the physics-informed neural network is good at capturing the dynamics of the actual solution. One caveat of this result is that the bound depends exponentially on $\|u\|_{C^0}$. Therefore, if the $C^0$-norm of the actual solution is very large, it can still be the case that even if the generalization error is small, the total error is not as small as expected. In this case, to achieve a small total error, the generalization error would need to be very small, which likely will incur additional computational costs. 
\begin{remark}
    Bounds for the constant of the Poincar\'e-Wirtinger are difficult to come by and depend on the geometry of the underlying domain. One prominent result shows that if $\Omega$ is a bounded convex Lipschitz domain with diameter $d$, then $C_{pw}\leq\frac{d}{\pi}$ \cite{Bebendorf}.
\end{remark}
\subsection{Question 3}
Lastly, question 3 pertains to the approximation of the total error in terms of the training error. In practice we only have access to the training error but not the generalization error. Based on the results from question 2, for a PINN $u_{\theta}$ we will derive an upper bound on sums of the total error and total error of the time derivative, in terms of the training error and and the number of quadrature points. Assuming that we find a global minimizer $u_{\theta^*}$ it will be shown that the total error of the global minimizer can be made arbitrarily small.

Here we use again the midpoint rule and the same set-up as in Corollary \ref{BountTraininError} with the training error being expressed in terms of the quadrature rule as in \ref{ErrorQuad}. We can now derive the following bound on the sum of the total error and the total error of the time-derivative in terms of the training error and the number of training points.
\begin{theorem}\label{Q3}
    Let $T>0, d\in\mathbb{N}$ and let $\Omega\subset\mathbb{R}^d$ be a Lipschitz domain such that there exists a exists an $N\in\mathbb{N}$ and a partition of $\Omega$ into $M\sim N^d$ cubes of side-length $\frac{1}{N}$ such that the midpoints of these cubes are inside $\Omega$. Let $u\in C^4(\Omega\times[0,T])$ be a classical solution of the semi-linear wave equation \ref{WaveBVP} with $a(x)\in C^2(\Omega)$ and (A1) and (A2) satisfied. Let $u_{\theta}$ be a PINN with parameters $\theta\in\Theta_{L,W,R}$. Let $\sigma$ denote a four times continuously differentiable activation function of the neural network and assume $\|\sigma\|_{C^n}\geq 1$. Then the following error bound holds:
    \begin{equation*}
    \begin{split}
        \int_{\Omega_T}|u_{\theta}-u|^2dxdt+\int_{\Omega_T}|u_{\theta_t}-u_t|^2dxdt&\leq\mathcal{C}(M)T\mathrm{exp}\left(T(\mathrm{max}\{1,2C_{pw}^2\}(1+\hat{C}+\frac{2\sqrt{T}}{C_{pw}^2}))\right)\\
        &=\mathcal{O}(\mathcal{E}_T(\theta,\mathcal{S})^2+M_t^{-\frac{2}{d}}+M_{PDE}^{-\frac{2}{d+1}}+M_s^{-\frac{1}{d}}),
    \end{split}        
    \end{equation*}
    where 
    \begin{equation*}
        \begin{split}
            \mathcal{C}(M)&=\mathrm{max}\{1,2C^2_{pw}\}\biggl( C_1M_t^{-\frac{2}{d}}+\mathcal{E}_T^{u_1}(\theta,\mathcal{S}_t)^2+C_2M_{PDE}^{-\frac{2}{d+1}}+\mathcal{E}_T^{PDE}(\theta,\mathcal{S}_{PDE})^2+C_3M_t^{-\frac{2}{d}}+\mathcal{E}_T^{\nabla u}(\theta,\mathcal{S}_t)^2\\
            &\quad+2\sqrt{T|\partial\Omega|}(16^L(d+1)^2(e^2W^3R\|\sigma\|_{C^1})^L+\|u\|_{C^1})(C_4M_s^{-\frac{1}{d}}+\mathcal{E}_T^{s,u_t}(\theta,\mathcal{S}_s))+C_5M_t^{-\frac{2}{d}}\\
            &\quad+\frac{2}{C^2_{pw}}\mathcal{E}_T^{u_0}(\theta,\mathcal{S}_t)^2\biggr)
        \end{split}
    \end{equation*}
    and $C_{pw}$ and $\hat{C}$ are the same constants as in Theorem \ref{Q3} and $C_1,C_2,C_3,C_4$ and $C_4$ are defined in the proof and depend on $d, L, W, R, \|\sigma\|_{C^n}$ and $\|u\|_{C^n}$.
\end{theorem}

The result of Theorem \ref{Q3} can be seen as an \textit{a posteriori} bound on the total error of a physics-informed neural network since it presupposes knowledge of the training error of the trained PINN. Moreover, when it is additionally assumed that the algorithm finds a global minimum of the optimization problem \ref{OptimizationProblem2} we can derive an \textit{a priori} result that indicates the global error can be made arbitrarily small in that case. The proof works by combining the answers to question 1, in particular Corollary \ref{BountTraininError} with Theorem \ref{Q3} and presupposes that the training set and hypothesis space is large enough.
\begin{theorem}\label{apriori}
    Let the assumptions of theorem \ref{Q3} be satisfied and assume additionally that there exists $\gamma\in[0,1)$ such that
    \begin{equation*}
        r=\frac{4\gamma}{d+1-4\gamma}.
    \end{equation*}
    Let $d\in\mathbb{N}$ and let $\Omega\subset\mathbb{R}^d$ be a Lipschitz-domain. Let $\hat{\Omega}$ denote a subset of $\mathbb{R}^d$ such that $\Omega\subset\subset\hat{\Omega}$ and $\hat{\Omega}=\Pi_{i=1}^d[a_i,b_i]$ for some $a_i,b_i\in\mathbb{Z},a_i<b_i$. Let $\varepsilon>0, T>0, k>2(18d+55)=:\eta$ and let $u\in H^{k+1}(\Omega\times(0,T))$ be a classical solution of the semi-linear wave equation \ref{WaveBVP}. Let the hypothesis space $\Theta$ satisfy 
    \[R\geq\varepsilon^{-\frac{1}{k-\eta}}\mathrm{ln}(\varepsilon^{-1})\], $W\geq C\varepsilon^{\frac{-(d+1)}{k-\eta}}$ and $L\geq 3, 24L-1\geq r$, where $C$ is specified in the proof. Let $u_{\theta^*(\mathcal{S})}$ be the PINN with tanh activation function that solves the optimization problem \ref{OptimizationProblem2} where the training set $\mathcal{S}$ satisfies $M_{PDE}\geq\varepsilon^{-\frac{-(d+1)(\eta-1)}{k-\eta}}, M_t\geq\varepsilon^{-\frac{-d(\eta-1)}{k-\eta}}$ and $M_s\geq\varepsilon^{-\frac{-2d(\eta-1)}{k-\eta}}$. Then it holds that
    \begin{equation*}
        \sqrt{\|u-u_{\theta^*(\mathcal{S})}\|^2_{L^2(\Omega\times[0,T])}+\|u_t-u_{\theta^*(\mathcal{S})_t}\|^2_{L^2(\Omega\times[0,T])}}=\mathcal{O}(\varepsilon).
    \end{equation*}
\end{theorem}

This concludes the section on theoretical results. All three questions that were laid out in the beginning were answered affirmatively. Theorem \ref{Q1} provides bounds on the residuals of a certain kind of tanh physics-informed neural network approximating the semilinear wave equation which were used to show that there exists a PINN such that the generalization and training error can be made arbitrarily small. Theorem \ref{Q2} provides an answer to question 2, that is, provided that the generalization error is small, it can be shown that the total error of a phyics-informed neural network approximating the semilinear wave equation is correspondingly small. Lastly, Theorem \ref{Q3} provided \textit{a posteriori} bounds on the total error in terms of the training error of the PINN, thereby answering question 3. Moreover, the \textit{a priori} bound in theorem \ref{apriori} demonstrates that under the assumption that the algorithm finds a global minimum, we can always construct a physics-informed neural network with a certain number of neurons and collocation points such that the total error is guaranteed to be arbitrarily small.

\section{Numerical Experiments}\label{NumericalExperiments}
In the previous section theoretical error bounds for physics-informed neural networks approximating semilinear wave equations were derived. In this section we seek to illustrate these bounds using numerical experiments to empirically test them on the following damped wave equation:
\begin{equation}\label{DampedWE}
    \begin{cases}
        u_{tt}-\Delta u +2\pi u_t=0 \qquad &\text{on }\Omega\times[0,0.5]\\
        u(x,y,0)=u_0(x,y)=\cos(\pi x)\cos(\pi y) &\text{on }\Omega\\
        u_t(x,y,0)=u_1(x,y)=0&\text{on }\Omega\\
        u(x,y,t)=0  &\text{on }\partial\Omega\times[0,0.5]
    \end{cases}
\end{equation}
with $\Omega=[-0.5,0.5]^2$. Note that here $a(x)=2\pi$.\\
It is easy to check that this problem is solved by the following function:
\begin{equation*}
    u(x,y,t)=e^{-\pi t}((\cos(\pi t) + \sin(\pi t))\cos(\pi x)\cos(\pi y).
\end{equation*}
This damped wave equation was chosen over the semilinear wave equation because it is easy to find an analytical solution for this problem in contrast to a semilinear wave problem, due to the presence of the nonlinearity in the latter. Having access to an analytical solution enables us to more accurately evaluate the total error in the $L^2$ norm without incurring further errors due to further numerical approximations. Moreover, the solution is clearly smooth enough to satisfy all the regularity conditions of the theorems derived in the previous section.\\
The bound we will be testing will be the one from Theorem \ref{Q3}, adapted to the damped wave equation \ref{DampedWE}. It differs slightly from the original as it does not include the constant $\hat{C}$ in the exponential due to the absence of the nonlinearity $f$. The bound then takes the following form:
\begin{equation*}
        \int_{\Omega_T}|u_{\theta}-u|^2dxdt+\int_{\Omega_T}|u_{\theta_t}-u_t|^2dxdt\leq\mathcal{C}(M)T\mathrm{exp}\left(T(\mathrm{max}\{1,2C_{pw}^2\}(1+\frac{2\sqrt{T}}{C_{pw}^2}))\right),      
\end{equation*}
where 
\begin{equation*}
    \begin{split}
        \mathcal{C}(M)&=\mathrm{max}\{1,2C^2_{pw}\}\biggl( C_{\Omega}\|\mathcal{R}_{u_1}^2\|_{C^2}M_t^{-\frac{2}{d}}+\mathcal{E}_T^{u_1}(\theta,\mathcal{S}_t)^2+C_{\Omega_T}\|\mathcal{R}_{PDE}^2\|_{C^2}M_{PDE}^{-\frac{2}{d+1}}+\mathcal{E}_T^{PDE}(\theta,\mathcal{S}_{PDE})^2\\
        &\quad+C_{\Omega}\|\mathcal{R}_{\nabla u}^2\|_{C^2}M_t^{-\frac{2}{d}}+\mathcal{E}_T^{\nabla u}(\theta,\mathcal{S}_t)^2+2\sqrt{T|\partial\Omega|}\|\nabla (u_{\theta}-u)\|_{C^0}\bigl(\sqrt{C_{\partial\Omega\times[0,T]}\|\mathcal{R}_{s,u_t}^2\|_{C^2}}M_s^{-\frac{1}{d}}\\
        &\quad+\mathcal{E}_T^{s,u_t}(\theta,\mathcal{S}_s)\bigr)+\frac{2C_{\Omega}}{C^2_{pw}}\|\mathcal{R}_{u_0}^2\|_{C^2}M_t^{-\frac{2}{d}}+\frac{2}{C^2_{pw}}\mathcal{E}_T^{u_0}(\theta,\mathcal{S}_t)^2\biggr).
    \end{split}
\end{equation*}

\begin{figure}
    \centering
    \subfloat{
    \includegraphics[width=.3\linewidth]{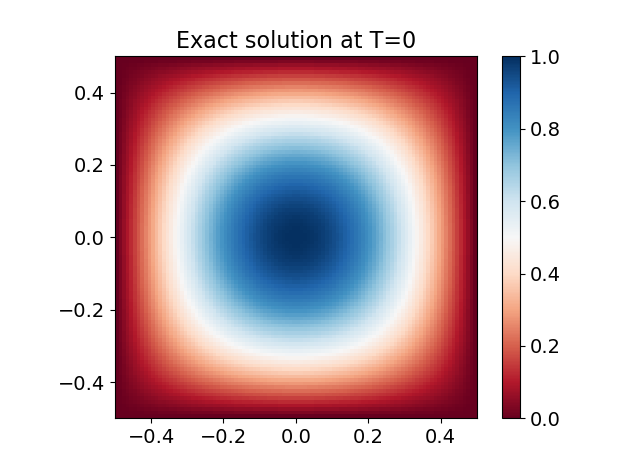}
    }
    \subfloat{
    \includegraphics[width=.3\linewidth]{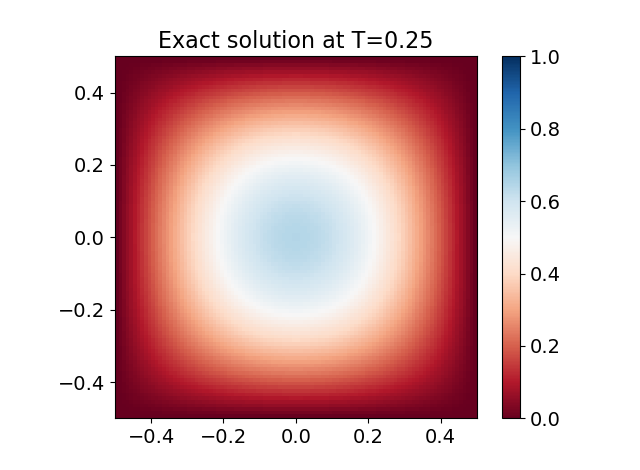}
    }
    \subfloat{
    \includegraphics[width=.3\linewidth]{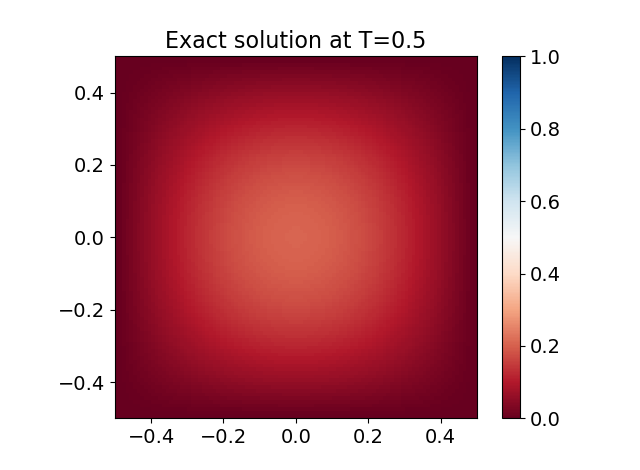}
    }\\
    \subfloat{
    \includegraphics[width=.3\linewidth]{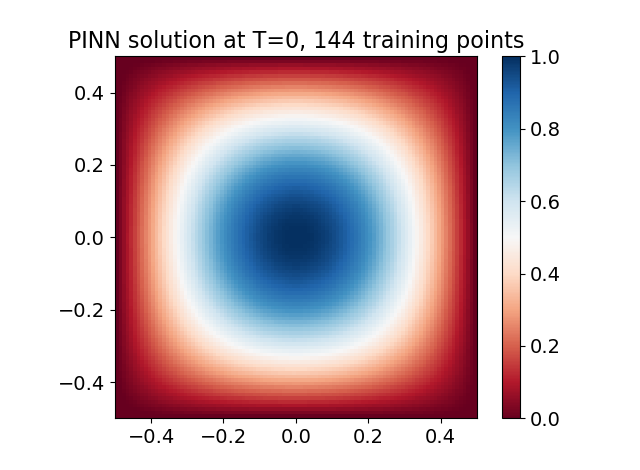}
    }
    \subfloat{
    \includegraphics[width=.3\linewidth]{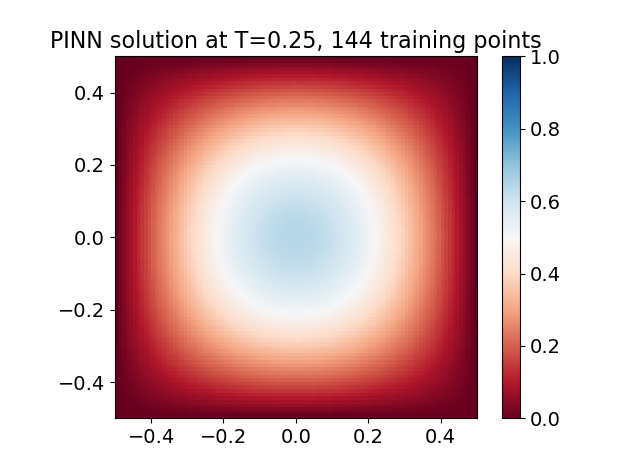}
    }
    \subfloat{
    \includegraphics[width=.3\linewidth]{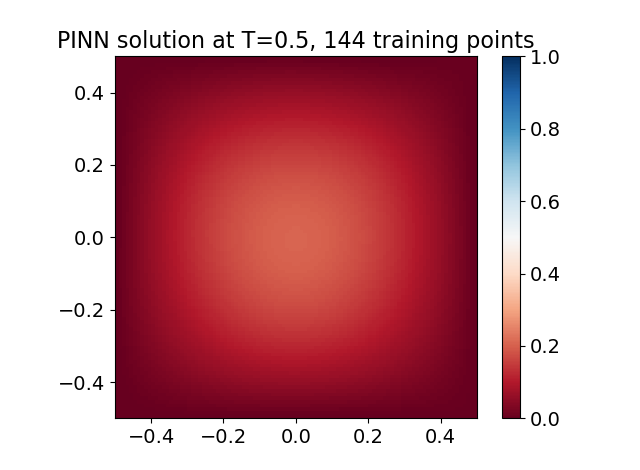}
    }\\
    \subfloat{
    \includegraphics[width=.3\linewidth]{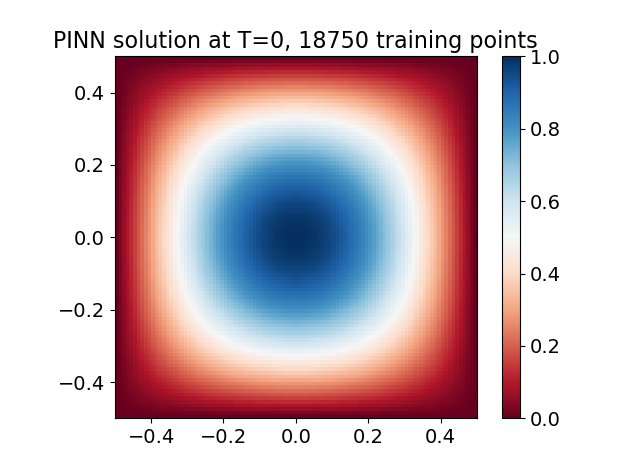}
    }
    \subfloat{
    \includegraphics[width=.3\linewidth]{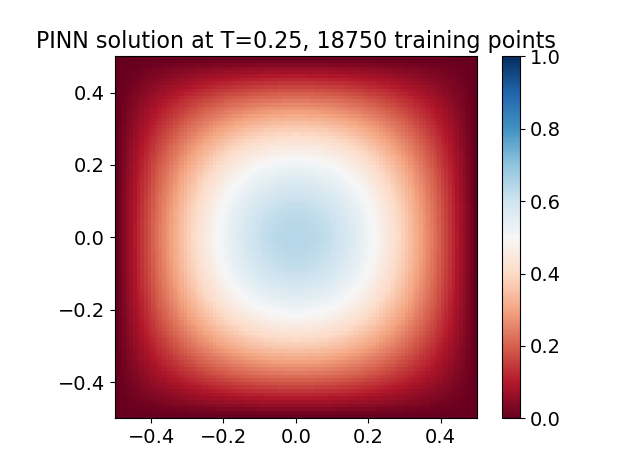}
    }
    \subfloat{
    \includegraphics[width=.3\linewidth]{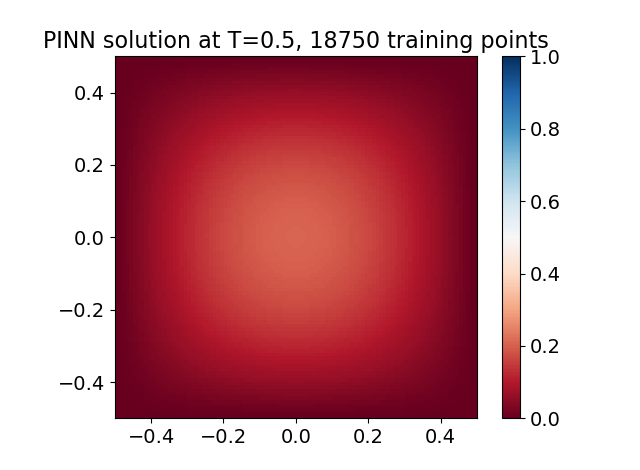}
    }
    
    \caption{First row: exact solution at T=0, T=0,25 and T=0.5. Second row: PINN solution with 144 training points at T=0, T=0.25 and T=0.5. Third row: PINN solution with 18750 training points at T=0, T=0.25 and T=0.5.}
    \label{fig:solutions}

    \subfloat{
    \includegraphics[width=.3\linewidth]{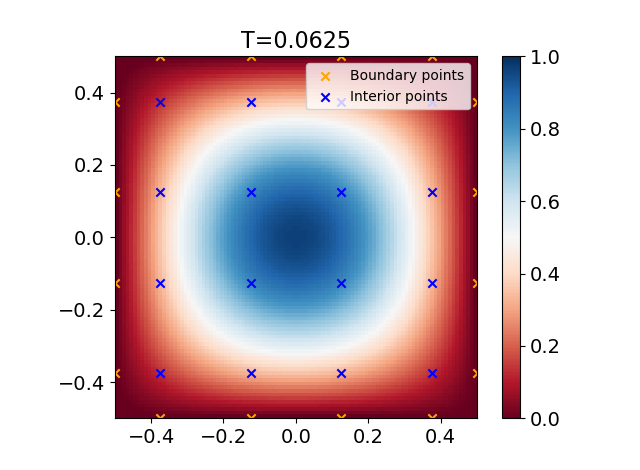}
    }
    \subfloat{
    \includegraphics[width=.3\linewidth]{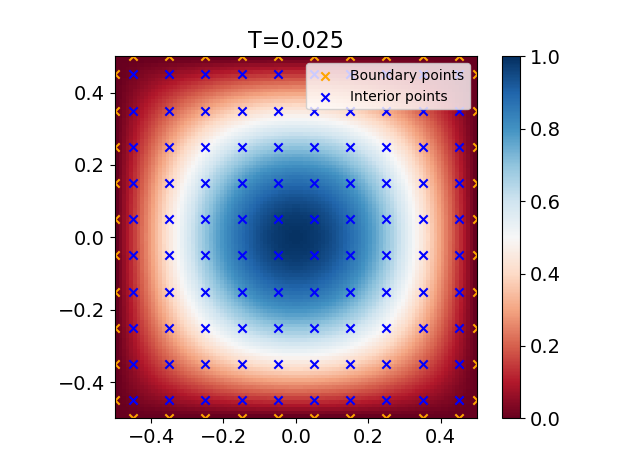}
    }
    \subfloat{
    \includegraphics[width=.3\linewidth]{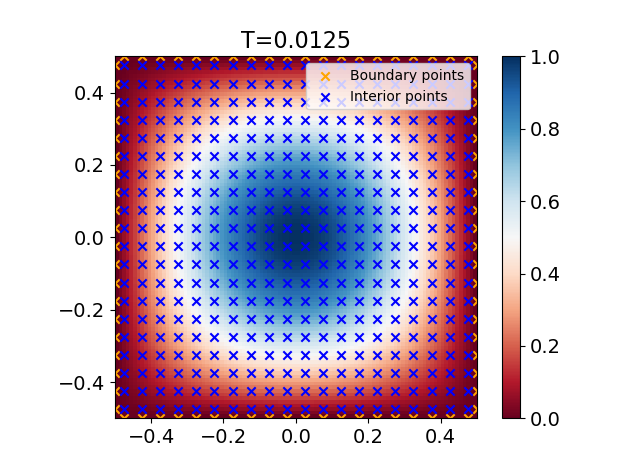}
    }
    
    \caption{Distribution of training points computed with the midpoint rule at the first time slice with PINN solution in the background for 144 total training points (left), 1500 total training points (middle) and 10000 total trainig points (right).}
    \label{fig:Trainingpoints}
\end{figure}

\begin{figure}
    \centering
    \subfloat{
    \includegraphics[width=.5\linewidth]{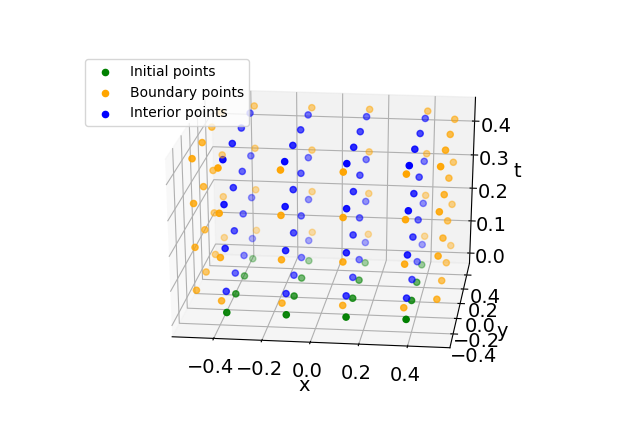}
    }
    \subfloat{
    \includegraphics[width=.5\linewidth]{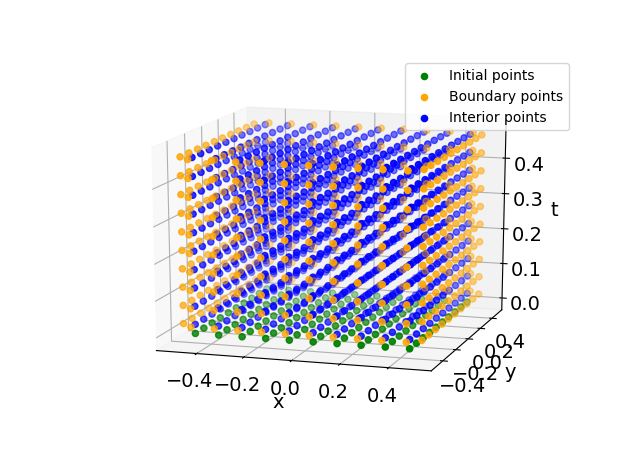}
    }
    \caption{Distribution of training points computed by using the midpoint rule in the space-time domain $[-0.5,0.5]\times[-0.5,0.5]\times[0,0.5]$ for 144 total training points (left) and 1500 total training points (right)}
    \label{fig:Trainingpoints3D}
    \centering
    \subfloat{
    \includegraphics[width=.3\linewidth]{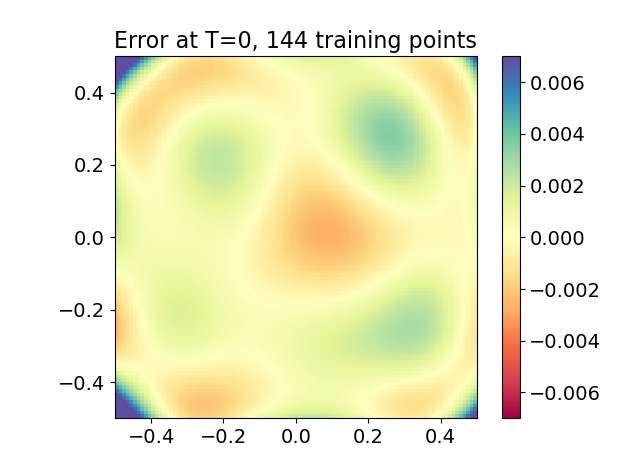}
    }
    \subfloat{
    \includegraphics[width=.3\linewidth]{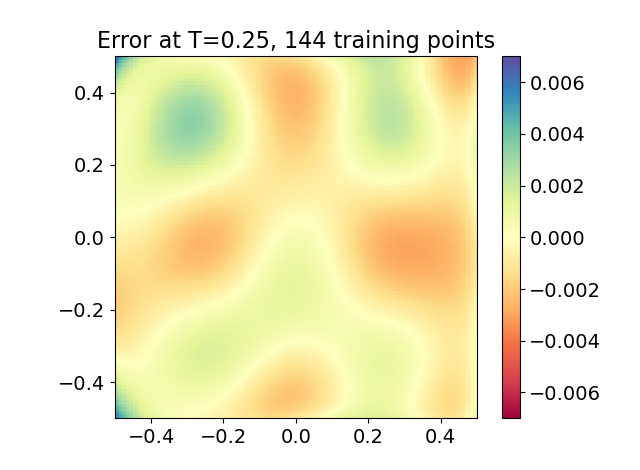}
    }
    \subfloat{
    \includegraphics[width=.3\linewidth]{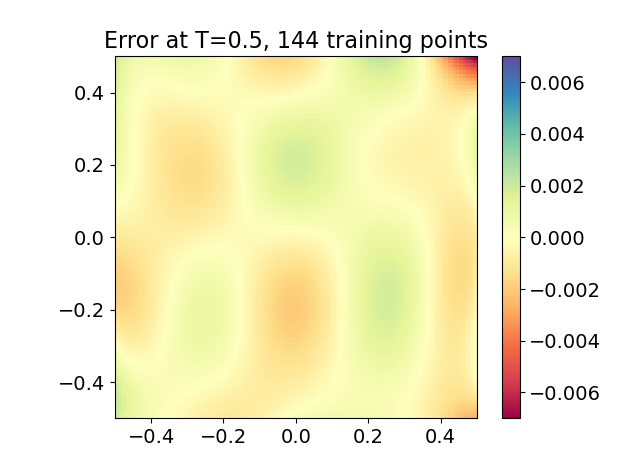}
    }\\
    \subfloat{
    \includegraphics[width=.3\linewidth]{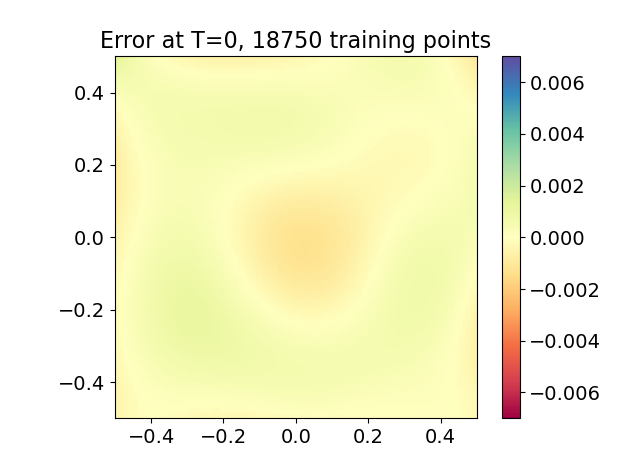}
    }
    \subfloat{
    \includegraphics[width=.3\linewidth]{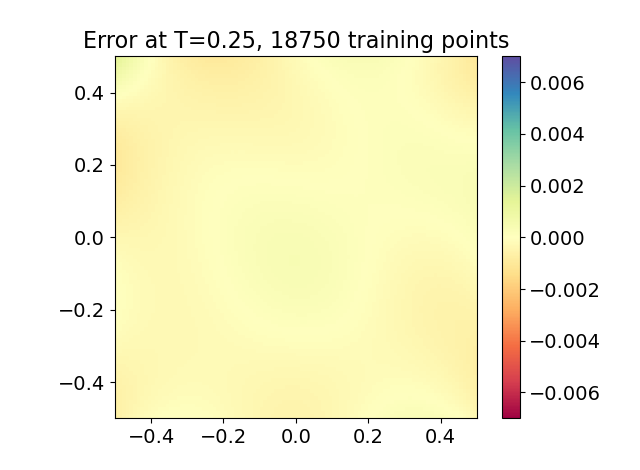}
    }
    \subfloat{
    \includegraphics[width=.3\linewidth]{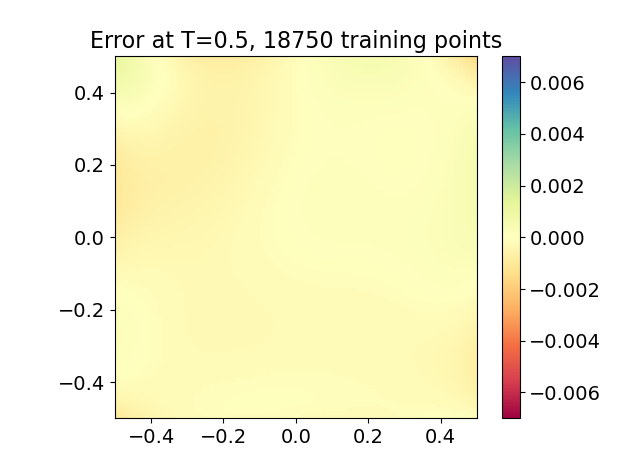}
    }
    
    \caption{First row: Error between exact and PINN solution at T=0, T=0.25 and T=0.5 for PINN with 144 training points. Second row: Error between exact and PINN solution at T=0, T=0.25 and T=0.5 for PINN with 18750 training points.}
    \label{fig:error}
    
    \centering
    \subfloat{
    \includegraphics[width=.3\linewidth]{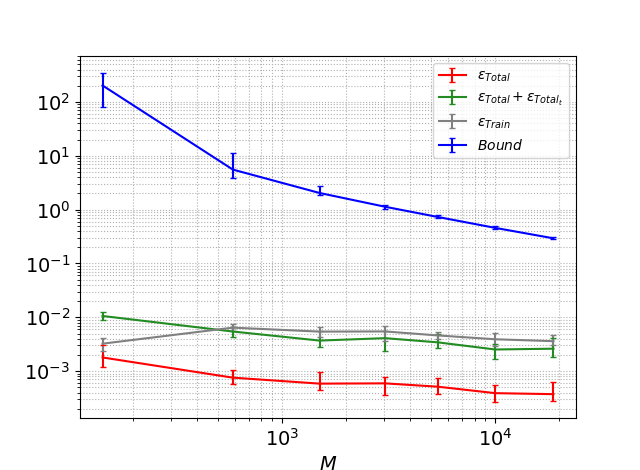}
    }   
    \subfloat{
    \includegraphics[width=.3\linewidth]{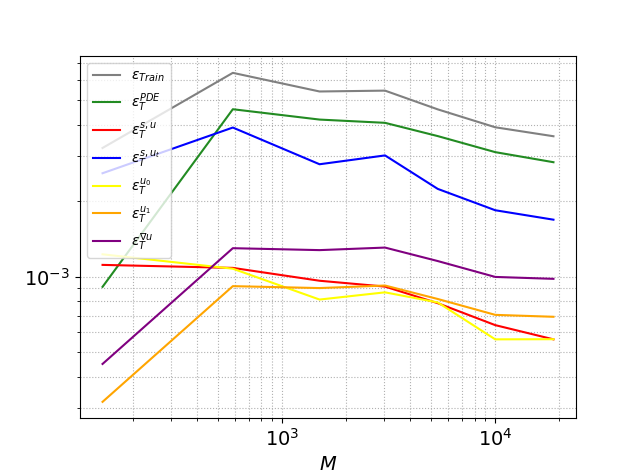}
    }
    \subfloat{
    \includegraphics[width=.3\linewidth]{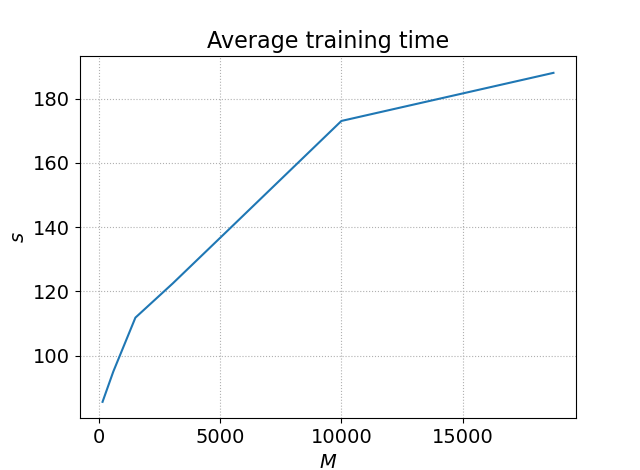}
    }
    
    \caption{On the left, we see the total error, training error and bound for different number of training points. In the center, we see the average training time per PINN in seconds over the number of training points. And on the right, we see the total training error as well as individual training errors over the number of training points.}
    \label{fig:Trainingerrors/time}
\end{figure}

The neural network used in the numerical experiments has two hidden layers and uses $tanh$ activation function (as suggested by the theory) with 80 neurons in each layer (for a graphical representation see figure \ref{fig:NN}). Full batch training with an L-BFGS optimizer for a maximum of 50000 iterations was used. The network was trained using empirical risk minimization on the training error \ref{TrainingErr} with the weights and training set determined by the midpoint rule. The grid-like distribution of the training points can be observed in figure \ref{fig:Trainingpoints} and \ref{fig:Trainingpoints3D}. For each number of training points the algorithm was run 10 times and the mean was computed. The code was adapted from the code used for the numerical experiments in \cite{Mishra2} by Mishra and Molinaro.

The exact solution at times $T=0, T=0.25$ and $T=0.5$, as well as the PINN solutions for the least and most amount of training points is shown in figure \ref{fig:solutions}. It is evident that the PINN solutions approximate the real solution very well as they are almost indistinguishable. Pointwise error plots are offered in figure \ref{fig:error} where it can be seen that the errors decrease with the number of training points. For 144 training points the magnitude of the pointwise error is between $10^{-3}$ and $10^{-2}$ whereas for 18750 training points the errors  are more around the magnitude $10^{-3}$.

The main result of the numerical experiments, i.e. the comparison between total error and the computed bound, is shown on the left of figure \ref{fig:Trainingerrors/time}. $M$ denotes the total number of training points with the leftmost points corresponding to 144 training points and the rightmost points corresponding to 18750 training points. It can be seen that the total error decreases with the number of training points until it arrives at around $3*10^{-4}$ for 18750 training points which suggests that the neural network approximates the actual solution rather well. The sum of the total error and the total error of the time derivative mimics the behaviour of the total error but is about one order of magnitude larger than the total error alone throughout. The bound is monotonically decreasing in the number of training points as suggested by the theory but is not quantitatively sharp. In the end, there are about three orders of magnitude difference between the total error and the bound and two orders of magnitude difference between the sum of the total error plus the total error of the time-derivative and the bound. To put these results into context, we note that the authors in \cite{DeRyck1} found differences of around two orders of magnitude between the total error and their bound for the Navier-Stokes equation. Considering that they did not have the additional difficulty of dealing with the second time-derivative, three orders of magnitude, respectively two for the sum of the total errors seems quite acceptable.

A more detailed analysis of the training error is offered in the middle plot of figure \ref{fig:Trainingerrors/time}. Here the training error is broken down into its individual components as specified in equation \ref{TrainingErr}. It can be noted that the biggest contributors to the overall training error are the error on the PDE dynamics in the interior ($\mathcal{E}_T^{PDE}$) as well as the error of the time derivative on the boundary ($\mathcal{E}_T^{s,u_t})$. Moreover, the rightmost plot in figure \ref{fig:Trainingerrors/time} shows the average training time over the number of training points which roughly doubles from least to most number of training points.



\section{Conclusion}\label{Conclusion}
This article has demonstrated that it is possible to derive theoretical error bounds for physics-informed neural networks approximating the semilinear wave equation and that these bounds hold true in practice while not being quantitatively sharp. We have used the approach developed in \cite{DeRyck1} to answer three questions that jointly imply that there exists a neural network such that the total error between the actual and the PINN solution can be made arbitrarily small.

Question 1, whether there exist neural networks with arbitrarily small generalization and training errors, was answered affirmatively in Theorem \ref{Q1} and the subsequent corollaries for a \textit{tanh} neural network with two hidden layers. Question 2, whether the total error is correspondingly small if the generalization error is small, was answered in Theorem \ref{Q2}. Moreover, we did not only find a bound on the total error, but the sum of the total error and the total error of the time derivative, indicating that the PINN solution also captures the temporal dynamics of the solution well. Furthermore, the question whether the total error is small given the training error is small and the training set is sufficiently large was answered positively in Theorem \ref{Q3}. This has given us an \textit{a posteriori} bound on the total error in terms of the training error and the number of training points. In addition to the \textit{a posteriori} bound, an \textit{a priori} bound on the total error was found in Theorem \ref{apriori}. The \textit{a posteriori} bound allows us to compute an explicit error bound for the physics-informed neural network after training is complete, whereas, under the assumption that the algorithm converges to the global minimum, the \textit{a priori} bound guarantees the existence of a \textit{tanh} neural network with a certain number of neurons and training points that has an arbitrarily small total error. However, the assumption that the algorithm finds the global minimum is very strong and can usually not be guaranteed in actuality. Lastly, the theoretical error bounds were tested empirically with a damped wave equation in Section \ref{NumericalExperiments}. It was found that the difference between the total error and the bound amounted to about three orders of magnitude which means that the bound is not quantitatively sharp. However, this was to be expected considering the results in \cite{DeRyck1} for the Navier-Stokes equation which did not even involve estimates for the second time-derivative.


While we think that this paper has shown that the approach conceived in \cite{DeRyck1} to answer these three questions consecutively to arrive at a bound for the total error, is quite powerful and leads to good results, nevertheless, there are some limitations. Firstly, the neural network still suffers from the curse of dimensionality as its width growth exponentially in the dimension (see Theorem \ref{Q1}) and secondly, the numerical experiments have shown that the error bound, while acceptable, is not sharp in practice. Moreover, in the future we hope to extend the numerical experiments to an actual semilinear wave equation that includes a nonlinearity. Futhermore, one could also argue that it is quite menial work to answer these question individually for each prominent class of partial differential equations. Previous work, such as \cite{Mishra1} and \cite{Mishra2} offer a more general approach to bound the total error in terms of the generalization error and the generalization error in terms of the training error but do not concern themselves with question 1. Consequently, the approach followed here improves on these prior results and offers a more complete analysis.

\appendix

\section{ }\label{se:appendix}
Results from analysis:

\begin{lemma}\label{PWICorollary}
    Let $\Omega\subset\mathbb{R}^d$ be a bounded connected Lipschitz domain, $1\leq p<\infty$ and $u\in W^{1,p}(\Omega)$. Choose $q$ such that $\frac{1}{q}+\frac{1}{p}=1$. Using the Poincar\'e-Wirtinger inequality, see, e.g., \cite{evans}, we obtain
    \begin{equation*}
        \begin{split}
            \|u\|_{L^p(\Omega)}&\leq C\|\nabla u\|_{L^p(\Omega)}+\|u_{\Omega}\|_{L^p(\Omega)}\\
            &=C\|\nabla u\|_{L^p(\Omega)}+|\Omega|^{\frac{1}{p}}|u_{\Omega}|\\
            &=C\|\nabla u\|_{L^p(\Omega)}+\frac{1}{|\Omega|^{\frac{1}{q}}}\Bigl|\int_{\Omega}u(x,t)dx\Bigr|\\
            &=C\|\nabla u\|_{L^p(\Omega)}+\frac{1}{|\Omega|^{\frac{1}{q}}}\Bigl|\int_{\Omega}\left(u(x,0)+\int_0^tu_t(x,r)dr\right)dx\Bigr|\\
        \end{split}
    \end{equation*}
    For $p=q=2$ we obtain further with H\"older
    \begin{equation*}
         \|u\|_{L^p(\Omega)}\leq C\|\nabla u\|_{L^p(\Omega)}+\left(\int_{\Omega}|u(x,0)|^2dx\right)^{\frac{1}{2}}+\sqrt{T}\left(\int_{\Omega}\int_0^t|u_t(x,r)|^2drdx\right)^{\frac{1}{2}}.
    \end{equation*}
\end{lemma}
\begin{lemma}[\cite{DeRyck2}]\label{A6}
    Let $d\in\mathbb{R},k\in\mathbb{N}_0,\Omega\subset\mathbb{R}^d$ and $f,g\in C^k(\Omega)$. Then it holds that
    \begin{equation*}
        \|fg\|_{C^k(\Omega)}\leq 2^k\|f\|_{C^k(\Omega)}\|g\|_{C^k(\Omega)}.
    \end{equation*}
\end{lemma}

\section{Appendix}
Results about neural networks:
\begin{theorem}[\cite{DeRyck1}]\label{tanh}
    Let $d,n\geq 2, m\geq 3, \delta>0, a_i,b_i\in\mathbb{Z}$ with $a_i<b_i$ for $1\leq i \leq d,\Omega=\Pi_{i=1}^d[a_i,b_i]$ and $f\in H^m(\Omega)$. Note that $P_{n,d}=\{\alpha\in\mathbb{N}_0^d:|\alpha|=n\}$ and $|P_{n,d}|=\binom{n+d-1}{n}$. Then, for every $N\in\mathbb{N}$ with $N>5$, there exists a tanh neural network $\hat{f}^N$ with two hidden layers, one of width at most $3\lceil\frac{m+n-2}{2}\rceil|P_{m-1,d+1}|+\sum_{i=1}^d(b_i-a_i)(N-1)$ and another of width at most $3\lceil\frac{d+n}{2}\rceil|P_{d+1,d+1}|N^d\Pi_{i=1}^d(b_i-a_i)$, such that for $l\in\{0,1,2\}$ it holds that 
    \begin{equation*}
        \|f-\hat{f}^N\|_{H^l(\Omega)}\leq 2^l3^dC_{l,m,d,f}(1+\delta)\mathrm{ln}^l\bigl(\beta_{l,\delta,d,f}N^{d+m+2}\bigr)N^{-m+l}
    \end{equation*}
    and where we define
    \begin{equation*}
        \begin{split}
            &\beta_{l,\delta,d,f}=\frac{2^{ld}5\mathrm{max}\{\Pi_{i=1}^d(b_i-a_i),d\}\mathrm{max}\{\|f\|_{W^{l,\infty}(\Omega)},1\}}{3^d\delta\mathrm{min}\{1,C_{l,m,d,f}\}}\\
            &C_{l,m,d,f}=\underset{0\leq i\leq l}{\mathrm{max}}\binom{d+i-1}{i}^{\frac{1}{2}}\frac{((m-i)!)^{\frac{1}{2}}}{(\lceil\frac{m-i}{d}\rceil!)^{\frac{d}{2}}}\left(\frac{3\sqrt{d}}{\pi}\right)^{m-i}|f|_{H^m}.
        \end{split}
    \end{equation*}
    Moreover, the weights $\hat{f}^N$ scale as $\mathcal{O}(N\mathrm{ln}(N)+N^{\kappa})$ with $\kappa=\mathrm{max}\{m^2,d(2+m+d)\}/n$.
\end{theorem}
\begin{lemma}[\cite{DeRyck1}]\label{CnNorm}
    Let $d,n,L,W,l_L\in\mathbb{N}$ and let $u_{\theta}:\mathbb{R}^{d+1}\rightarrow\mathbb{R}^{l_L}$ be a neural network with $\theta\in\Theta_{L,W,R}$ for $L\geq 2, R,W\geq 1$. Assume that for the activation function $\sigma$ of the neural network it holds that $\|\sigma\|_{C^n}\geq 1$. Then it holds for $1\leq j\leq l_L$ that
    \begin{equation*}
        \|(u_{\theta})_j\|_{C^n}\leq 16^L(d+1)^{2n}\bigl(e^2n^4W^3R^n\|\sigma\|_{C^n}\bigr)^{nL}.
    \end{equation*}
\end{lemma}

\begin{lemma}\label{C0Norm}
     Let $d,L,W,l_L\in\mathbb{N}$ and let $u_{\theta}:\mathbb{R}^{d+1}\rightarrow\mathbb{R}^{l_L}$ be a neural network with $\theta\in\Theta_{L,W,R}$ for $L\geq 2, R,W\geq 1$. Assume that for the activation function $\sigma$ of the neural network it holds that $\|\sigma\|_{C^0}\geq 1$. Then it holds for $1\leq j\leq l_L$ that
     \begin{equation*}
         \|(u_{\theta})_j\|_{C^0}\leq R(W\|\sigma\|_{C^0}+1).
     \end{equation*}
\end{lemma}
\begin{proof}
    From the definition of a neural network \ref{NeuralNetwork} it follows that
    \begin{equation*}
    \begin{split}
        \|u_{\theta}\|_{C^0}&=\|f_L^{\theta}\circ f_{L-1}^{\theta}\circ...\circ f_1^{\theta}\|_{C^0}=\|W_L(f_{L-1}^{\theta}\circ...\circ f_1^{\theta})+b_L\|_{C^0}\\
        &=\|W_L(\sigma(W_{L-1}(f_{L-2}^{\theta}\circ...\circ f_1^{\theta})+b_{L-1}))+b_L\|_{C^0}\\
        &\leq R(W\|\sigma\|_{C^0}+1).
    \end{split}
    \end{equation*}
\end{proof}

\section{Proofs of the error bounds}
In this section, we present the proofs of the main results.

\begin{proof}[Proof of Theorem \ref{Q1}]
The first statement follows directly from \ref{Regularity}.

For the second statement, we first note that since $\Omega\subset\subset\hat{\Omega}$ and $\Omega$ Lipschitz, also $\Omega\times(0,T)\subset\subset\hat{\Omega}\times(-\varepsilon,T+\varepsilon)$ for some $\varepsilon>0$. We can use the Calderon Extension Theorem, see, e.g., \cite{Dobrowolski}, to extend any solution $u\in H^{k+1}(\Omega\times(0,T))$ to a function $\tilde{u}\in H^{k+1}_0(\hat{\Omega}\times(-\varepsilon,T+\varepsilon))$ with $\tilde{u}|_{\Omega\times(0,T)}=u$ and zero trace on the boundary $\partial(\hat{\Omega}\times(-\varepsilon,T+\varepsilon))$. We can then restrict $\tilde{u}$ again to $\hat{u}=\tilde{u}|_{\hat{\Omega}\times(0,T)}$ and thus $\hat{u}\in H^{k+1}(\hat{\Omega}\times(0,T))$ with $\hat{u}|_{\Omega\times(0,T)}=u$. We will continue to write $u$ for $\hat{u}$ when the context is clear.

By Theorem \ref{tanh}, for every $N>5$ there exists a tanh neural network $u_{\theta}$ with two hidden layers of width $3\lceil \frac{k+n-1}{2}\rceil\binom{d+k+1}{d}+\lceil (N-1)(T+\sum_{i=1}^d(b_i-a_i) \rceil$ and $3(d+3)N^{d+1}\lceil \frac{d+n+1}{2} \rceil \lceil T\Pi_{i=1}^d(b_i-a_i) \rceil$ such that for every $0\leq l\leq 2$
\begin{equation*}
    \|u-u_{\theta}\|_{H^l(\hat{\Omega}_T)}\leq C_l\lambda_l(N)N^{-k-1+l}
\end{equation*}
and the weights scale as $\mathcal{O}(N\mathrm{ln}(N)+N^{\kappa})$, where $\kappa=\mathrm{max}\{(k+1)^2,(d+1)(d+k+4)\}/n$.\\
We can use this result to bound the residuals, where we start with the PDE residual $\mathcal{R}_{PDE}$.\\
We have
\begin{equation*}
    \|u_{tt}-u_{\theta_{tt}}\|_{L^2(\Omega_T)}\leq\|u_{tt}-u_{\theta_{tt}}\|_{L^2(\hat{\Omega}_T)}\leq\|u-u_{\theta}\|_{H^2(\hat{\Omega}_T)}\leq C_2\lambda_2(N)N^{-k+1}.
\end{equation*}
Moreover,
\begin{equation*}
     \|\Delta u-\Delta u_{\theta}\|_{L^2(\Omega_T)}\leq\|\Delta u-\Delta u_{\theta}\|_{L^2(\hat{\Omega}_T)}
    \leq\sqrt{d}\|u-u_{\theta}\|_{H^2(\hat{\Omega}_T)}\leq \sqrt{d}C_2\lambda_2(N)N^{-k+1},
\end{equation*}
and
\begin{equation*}
\begin{split}
    &\|a(x)(u_t-u_{\theta_t})\|_{L^2(\Omega_T)}\leq\|a\|_{L^{\infty}(\Omega)}\|u_t-u_{\theta_t}\|_{L^2(\Omega_T)}\leq\|a\|_{L^{\infty}(\Omega)}\|u_t-u_{\theta_t}\|_{L^2(\hat{\Omega}_T)}\\
    &\leq \|a\|_{L^{\infty}(\Omega)}\|u-u_{\theta}\|_{H^1(\hat{\Omega}_T)}\leq \|a\|_{L^{\infty}(\Omega)}C_1\lambda_1(N)N^{-k}.
\end{split}
\end{equation*}
For $f(x,u)$ we use the Mean Value Theorem for G\^ateaux Differentiable Functions, see, e.g., \cite{MVT}, and assumption (A2) to get
\begin{align}\label{fEstimation}
        &\|f(x,u)-f(x,u_{\theta})\|_{L^2(\Omega_T)}^2 \notag \\
        &=\int_0^T\int_{\Omega}|f(x,u)-f(x,u_{\theta})|^2dxdt \notag\\
        &=\int_0^T\int_{\Omega}\left|\int_0^1\frac{\partial}{\partial u}f(x,u-s(u_{\theta}-u))(u_{\theta}-u)ds\right|^2dxdt \notag\\
        &\leq\int_0^T\int_{\Omega}\left|\int_0^1c|u-s(u_{\theta}-u)|^r|u_{\theta}-u|ds\right|^2dxdt \notag\\
        &\leq\int_0^T\int_{\Omega}\left|c2^{r-1}|u_{\theta}-u|\int_0^1|u|^r+s^r|u_{\theta}-u|^rds\right|^2dxdt \notag\\
        &\leq\int_0^T\int_{\Omega}\left|c2^{r-1}(|u|^r+\frac{1}{r+1}|u_{\theta}-u|^r)|u_{\theta}-u|\right|^2dxdt \notag\\
        &\leq\int_0^T\int_{\Omega}c^22^{2r-1}(|u|^{2r}+\frac{1}{(r+1)^2}|u_{\theta}-u|^{2r})|u_{\theta}-u|^2dxdt \notag\\
        &\leq c^22^{2r-1}\left(\int_0^T\int_{\Omega}(|u|^{2r}+\frac{1}{(r+1)^2}|u_{\theta}-u|^{2r})^pdxdt\right)^{\frac{1}{p}}\left(\int_0^T\int_{\Omega}|u_{\theta}-u|^{2q}dxdt\right)^{\frac{1}{q}} \notag\\
        &\leq c^22^{2r-\frac{1}{p}}\left(\int_0^T\int_{\Omega}(|u|^{2pr}+\frac{1}{(r+1)^{2p}}|u_{\theta}-u|^{2pr})dxdt\right)^{\frac{1}{p}}\left(\int_0^T\int_{\Omega}|u_{\theta}-u|^{2q}dxdt\right)^{\frac{1}{q}},
\end{align}
where in the second to last line we used the H\"older inequality. We choose $p$ and $q$ such that $2pr=2q$, i.e., $q=pr=\frac{q}{q-1}r$ implying $q=r+1$. Using the Gagliardo-Nirenberg inequality, see, e.g., \cite{Brezis},  with $\gamma=\frac{r(d+1)}{4(r+1)}$ we then get
\begin{equation*}
    \begin{split}
    &= c^22^{\frac{2r^2+r}{r+1}}\left(\int_0^T\int_{\Omega}(|u|^{2(r+1)}+\frac{1}{(r+1)^{2+\frac{2}{r}}}|u_{\theta}-u|^{2(r+1)})dxdt\right)^{\frac{r}{r+1}}\left(\int_0^T\int_{\Omega}|u_{\theta}-u|^{2(r+1)}dxdt\right)^{\frac{1}{r+1}}\\
    &\leq c^22^{2r-1}\left(\|u\|_{L^{2(r+1)}(\Omega_T)}^{2r}+\frac{1}{(r+1)^2}\|u_{\theta}-u\|_{L^{2(r+1)}(\Omega_T)}^{2r}\right)\|u_{\theta}-u\|_{L^{2(r+1)}(\Omega_T)}^{2}\\
    &\leq c^22^{2r-1}C^2\left(\|u\|_{H^2(\Omega_T)}^{{2r\gamma}}\|u\|_{L^2(\Omega_T)}^{{2r(1-\gamma)}}+\frac{1}{(r+1)^2}\|u_{\theta}-u\|_{H^2(\Omega_T)}^{{2r\gamma}}\|u_{\theta}-u\|_{L^2(\Omega_T)}^{{2r(1-\gamma)}}\right)\\
    &\quad\|u_{\theta}-u\|_{H^2(\Omega_T)}^{{2\gamma}}\|u_{\theta}-u\|_{L^2(\Omega_T)}^{{2(1-\gamma)}}\\
    &\leq c^22^{2r-1}C^2\left(\|u\|_{H^2(\Omega_T)}^{{2r\gamma}}\|u\|_{L^2(\Omega_T)}^{{2r(1-\gamma)}}+\frac{1}{(r+1)^2}\|u_{\theta}-u\|_{H^2(\hat{\Omega}_T)}^{{2r\gamma}}\|u_{\theta}-u\|_{L^2(\hat{\Omega}_T)}^{{2r(1-\gamma)}}\right)\\
    &\quad\|u_{\theta}-u\|_{H^2(\hat{\Omega}_T)}^{{2\gamma}}\|u_{\theta}-u\|_{L^2(\hat{\Omega}_T)}^{{2(1-\gamma)}}\\
    &\leq c^22^{2r-1}C^2\left(\|u\|_{H^2(\Omega_T)}^{{2r\gamma}}\|u\|_{L^2(\Omega_T)}^{{2r(1-\gamma)}}+\frac{1}{(r+1)^2}(C_2\lambda_2(N)N^{-k+1})^{{2r\gamma}}(C_0\lambda_0(N)N^{-k-1})^{{2r(1-\gamma)}}\right)\\
    &\quad(C_2\lambda_2(N)N^{-k+1})^{{2\gamma}}(C_0\lambda_0(N)N^{-k-1})^{{2(1-\gamma)}},
    \end{split}
\end{equation*}
where $C$ is the constant that we get from the Gagliardo-Nirenberg inequality. We then arrive at 

\begin{align*}
\|f(x,u)-f(x,u_{\theta})\|_{L^2(\Omega_T)}
&\leq c\,2^{r-\frac{1}{2}}C
\Bigl(
\|u\|_{H^2(\Omega_T)}^{2r\gamma}
\|u\|_{L^2(\Omega_T)}^{2r(1-\gamma)}
\\
&\qquad
+\frac{1}{(r+1)^2}
(C_2\lambda_2(N)N^{-k+1})^{2r\gamma}
(C_0\lambda_0(N)N^{-k-1})^{2r(1-\gamma)}
\\
&\qquad
(C_2\lambda_2(N)N^{-k+1})^{\gamma}
(C_0\lambda_0(N)N^{-k-1})^{1-\gamma}
\Bigr)^{1/2}
\\
&=:\tilde{C}
(C_2\lambda_2(N)N^{-k+1})^{\gamma}
(C_0\lambda_0(N)N^{-k-1})^{1-\gamma}.
\end{align*}

All in all, using the triangle inequality and that $\|\mathcal{R}_{PDE}[u]\|_{L^2(\Omega_T)}=0$ for the exact solution $u$, we get
\begin{equation*}
    \begin{split}
        &\|\mathcal{R}_{PDE}[u_{\theta}]\|_{L^2(\Omega_T)}=\|u_{\theta_{tt}}-\Delta u_{\theta}+a(x)u_{\theta_t}+f(x,u_{\theta})\|_{L^2(\Omega_T)}\\
        &=\|u_{\theta_{tt}}-u_{tt}-\Delta u_{\theta}+\Delta u+a(x)u_{\theta_t}-a(x)u+f(x,u_{\theta})-f(x,u)+u_{tt}-\Delta u+a(x)u+f(x,u)\|_{L^2(\Omega_T)}\\
        &\leq\|u_{\theta_{tt}}-u_{tt}\|_{L^2(\Omega_T)}+\|\Delta u-\Delta u_{\theta}\|_{L^2(\Omega_T)}+\|a(x)u_{\theta_t}-a(x)u\|_{L^2(\Omega_T)}+\|f(x,u_{\theta})-f(x,u)\|_{L^2(\Omega_T)}\\
        &\leq (C_2+\sqrt{d})\lambda_2(N)N^{-k+1}+\|a\|_{L^{\infty}(\Omega)}C_1\lambda_1(N)N^{-k}\\
        &\quad+\tilde{C}(C_2\lambda_2(N)N^{-k+1})^{\frac{r(d+1)}{4(r+1)}}(C_0\lambda_0(N)N^{-k-1})^{\frac{4(r+1)-r(d+1)}{4(r+1)}}.
    \end{split}
\end{equation*}
Similarly, we use the multiplicative trace inequality, see \cite{DeRyck1}, to find an upper bound for the boundary residuals.
\begin{equation*}
    \begin{split}
        &\|\mathcal{R}_{s,u}[u_{\theta}]\|_{L^2(\partial\Omega\times[0,T])}=
        \|u_{\theta}\|_{L^2(\partial\Omega\times[0,T])}\leq\|u-u_{\theta}\|_{L^2(\partial\Omega\times[0,T])}+\|u\|_{L^2(\partial\Omega\times[0,T])}\\
        &=\|u-u_{\theta}\|_{L^2(\partial\Omega\times[0,T])}\leq\|u-u_{\theta}\|_{L^2(\partial\Omega_T)}\leq\sqrt{\frac{2\mathrm{max}\{2h_{\Omega_T},d+1\}}{\rho_{\Omega_T}}}\|u-u_{\theta}\|_{H^1(\Omega_T)}\\
        &\leq\sqrt{\frac{2\mathrm{max}\{2h_{\Omega_T},d+1\}}{\rho_{\Omega_T}}}\|u-u_{\theta}\|_{H^1(\hat{\Omega}_T)}\leq\sqrt{\frac{2\mathrm{max}\{2h_{\Omega_T},d+1\}}{\rho_{\Omega_T}}}C_1\lambda_1(N)N^{-k},
    \end{split}
\end{equation*}
where $h_{\Omega_T}$ denotes the diameter of $\Omega_T$ and $\rho_{\Omega_T}$ the radius of the largest $d+1$-dimensional ball that can be inscribed into $\Omega_T$.

In the same way we find the bound
\begin{equation*}
    \begin{split}
        &\|\mathcal{R}_{s,u_t}[u_{\theta}]\|_{L^2(\partial\Omega\times[0,T])}\leq\sqrt{\frac{2\mathrm{max}\{2h_{\Omega_T},d+1\}}{\rho_{\Omega_T}}}\|u_t-u_{\theta_t}\|_{H^1(\Omega_T)}\\
        &\leq\sqrt{\frac{2\mathrm{max}\{2h_{\Omega_T},d+1\}}{\rho_{\Omega_T}}}\|u-u_{\theta}\|_{H^2(\hat{\Omega}_T)}\leq\sqrt{\frac{2\mathrm{max}\{2h_{\Omega_T},d+1\}}{\rho_{\Omega_T}}}C_2\lambda_2(N)N^{-k+1}.
    \end{split}
\end{equation*}
To bound the residuals associated with the initial values we proceed in a similar manner using the multiplicative trace inequality, see \cite{DeRyck1}, and the triangle inequality and get 
\begin{equation*}
    \begin{split}
        &\|\mathcal{R}_{u_0}[u_{\theta}]\|_{L^2(\Omega)}=\|u_0-u_{\theta}(t=0)\|_{L^2(\Omega)}\leq\|u(t=0)-u_{\theta}(t=0)\|_{L^2(\Omega)}+\|u_0-u(t=0)\|_{L^2(\Omega)}\\
        &=\|u(t=0)-u_{\theta}(t=0)\|_{L^2(\Omega)}\leq\|u-u_{\theta}\|_{L^2(\partial\Omega_T)}\leq\sqrt{\frac{2\text{max}\{2h_{\Omega_T},d+1\}}{\rho_{\Omega_T}}}\|u-u_{\theta}\|_{H^1(\Omega_T)}\\
        &\leq\sqrt{\frac{2\text{max}\{2h_{\Omega_T},d+1\}}{\rho_{\Omega_T}}}\|u-u_{\theta}\|_{H^1(\hat{\Omega}_T)}\leq\sqrt{\frac{2\text{max}\{2h_{\Omega_T},d+1\}}{\rho_{\Omega_T}}}C_1\lambda_1(N)N^{-k}
    \end{split}
\end{equation*}
and
\begin{equation*}
    \begin{split}
        &\|\mathcal{R}_{u_1}[u_{\theta}]\|_{L^2(\Omega)}=\|u_1-u_{\theta_t}(t=0)\|_{L^2(\Omega)}\leq\|u_t-u_{\theta_t}\|_{L^2(\partial\Omega_T)}\leq\sqrt{\frac{2\text{max}\{2h_{\Omega_T},d+1\}}{\rho_{\Omega_T}}}\|u_t-u_{\theta_t}\|_{H^1(\Omega_T)}\\
        &\leq\sqrt{\frac{2\text{max}\{2h_{\Omega_T},d+1\}}{\rho_{\Omega_T}}}\|u-u_{\theta}\|_{H^2(\hat{\Omega}_T)}\leq\sqrt{\frac{2\text{max}\{2h_{\Omega_T},d+1\}}{\rho_{\Omega_T}}}C_2\lambda_2(N)N^{-k+1}.
    \end{split}
\end{equation*}
Lastly, we find that
\begin{equation*}
    \begin{split}
        &\|\mathcal{R}_{\nabla u}[u_{\theta}]\|_{L^2(\Omega)}^2=\|\nabla u_0-\nabla u_{\theta}(t=0)\|_{L^2(\Omega)}^2\leq\|\nabla u (t=0)-\nabla u_{\theta}(t=0)\|_{L^2(\Omega)}^2\\
        &=\sum_{i=1}^d\|\frac{\partial}{\partial x_i}u(t=0)-\frac{\partial}{\partial x_i}u_{\theta}(t=0)\|_{L^2(\Omega)}^2\leq\sum_{i=1}^d\|\frac{\partial}{\partial x_i}u-\frac{\partial}{\partial x_i}u_{\theta}\|_{L^2(\partial\Omega_T)}^2\\
        &\leq\sum_{i=1}^d\frac{2\text{max}\{2h_{\Omega_T},d+1\}}{\rho_{\Omega_T}}\left\Vert\frac{\partial}{\partial x_i}u-\frac{\partial}{\partial x_i}u_{\theta}\right\Vert_{H^1(\Omega_T)}^2\leq\sum_{i=1}^d\frac{2\text{max}\{2h_{\Omega_T},d+1\}}{\rho_{\Omega_T}}\|u-u_{\theta}\|_{H^2(\Omega_T)}^2\\
        &=\frac{2d\text{max}\{2h_{\Omega_T},d+1\}}{\rho_{\Omega_T}}\|u-u_{\theta}\|_{H^2(\Omega_T)}^2\leq\frac{2d\text{max}\{2h_{\Omega_T},d+1\}}{\rho_{\Omega_T}}\|u-u_{\theta}\|_{H^2(\hat{\Omega}_T)}^2
    \end{split}
\end{equation*}
and thus
\begin{equation*}
    \|\mathcal{R}_{\nabla u}[u_{\theta}]\|_{L^2(\Omega)}\leq\sqrt{\frac{2d\text{max}\{2h_{\hat{\Omega}_T},d+1\}}{\rho_{\Omega_T}}}\|u-u_{\theta}\|_{H^2(\hat{\Omega}_T)}\leq\sqrt{\frac{2d\text{max}\{2h_{\Omega_T},d+1\}}{\rho_{\Omega_T}}}C_2\lambda_2(N)N^{-k+1}.
\end{equation*}
\end{proof}

\begin{proof}[Proof of Theorem \ref{Q2}]
    Let $\hat{u}=u_\theta-u$ denote the difference between the solution to the semi-linear wave equation and the PINN with parameter $\theta$. We can write the residuals as
\begin{equation*}
    \begin{split}
        &\mathcal{R}_{PDE}=\hat{u}_{tt}-\Delta\hat{u}+a(x)\hat{u}_t+f(x,u_{\theta})-f(x,u),\\
        &\mathcal{R}_{u_0}=\hat{u}(t=0),\\
        &\mathcal{R}_{u_1}=\hat{u_t}(t=0),\\
        &\mathcal{R}_{S,u}=\hat{u}(x,t),\\
        &\mathcal{R}_{S,u_t}=\hat{u}_t(x,t),\\
        &\mathcal{R}_{\nabla u}=\nabla\hat{u}(x,0).
    \end{split}
\end{equation*}
We note that $\frac{\partial}{\partial t}|\hat{u}_t|^2=2\hat{u}_{tt}\hat{u}_t$. Multiplying $\mathcal{R}_{PDE}$ by $2\hat{u}_t$ we get
\begin{equation*}
2\hat{u}_t\mathcal{R}_{PDE}=2\hat{u}_t\hat{u}_{tt}-2\hat{u}_t\Delta\hat{u}+2a(x)|\hat{u}_t|^2+2\hat{u}_t(f(x,u_{\theta})-f(x,u)).
\end{equation*}
Using assumption (A1) that $a(x)\geq0, a(x)\neq0$ a.e., we get
\begin{equation*}
    \begin{split}
        \frac{\partial}{\partial t}|\hat{u}_t|^2&=2\hat{u}_t\mathcal{R}_{PDE}+2\hat{u}_t\Delta\hat{u}-2a(x)|\hat{u}_t|^2-2\hat{u}_t(f(x,u_{\theta})-f(x,u))\\
        &\leq2\hat{u}_t\mathcal{R}_{PDE}+2\hat{u}_t\Delta\hat{u}-2\hat{u}_t(f(x,u_{\theta})-f(x,u))
    \end{split}
\end{equation*}
Using that
\begin{equation*}
\begin{split}
        \frac{d}{dt}\int_{\Omega}|\nabla\hat{u}|^2dx&=\int_{\Omega}\frac{\partial}{\partial t}|\nabla\hat{u}|^2dx=2\int_{\Omega}\nabla\hat{u}\cdot\nabla\hat{u}_tdx\\
        &=-2\int_{\Omega}\hat{u}_t\Delta\hat{u}dx+2\int_{\partial\Omega}\hat{u}_t(\hat{n}\cdot\nabla\hat{u})ds(x),
\end{split}
\end{equation*}
we integrate over $\Omega$ and have
\begin{equation}\label{Zwischenergebnis1}
    \begin{split}
        &\frac{d}{dt}\int_{\Omega}|\hat{u}_t|^2dx\leq\int_{\Omega}2\hat{u}_t\mathcal{R}_{PDE}+2\hat{u}_t\Delta\hat{u}-2\hat{u}_t(f(x,u_{\theta})-f(x,u))dx\\
        &=\int_{\Omega}2\hat{u}_t\mathcal{R}_{PDE}dx-\frac{d}{dt}\int_{\Omega}|\nabla\hat{u}|^2-2\hat{u}_t(f(x,u_{\theta})-f(x,u))dx+2\int_{\partial\Omega}\hat{u}_t(\hat{n}\cdot\nabla\hat{u})ds(x)\\
        &\leq\int_{\Omega}2\hat{u}_t\mathcal{R}_{PDE}dx-\frac{d}{dt}\int_{\Omega}|\nabla\hat{u}|^2+2|\hat{u}_t||f(x,u_{\theta})-f(x,u)|dx+2\int_{\partial\Omega}\hat{u}_t(\hat{n}\cdot\nabla\hat{u})ds(x).
    \end{split}
\end{equation}
Using Cauchy-Schwarz, the boundary term can be estimated as follows:
\begin{equation}\label{boundary}
    \begin{split}
        &\int_{\partial\Omega}\hat{u}_t(\hat{n}\cdot\nabla\hat{u})ds(x)\leq\int_{\partial\Omega}|\hat{u}_t||\hat{n}||\nabla\hat{u}|ds(x)\\
        &\leq\|\nabla\hat{u}\|_{C^0(\partial\Omega)}\int_{\partial\Omega}|\hat{u}_t|ds(x)\leq\|\hat{u}\|_{C^1}\|\mathcal{R}_{S,u_t}\|_{L^1(\partial\Omega)}
    \end{split}
\end{equation}
The term $2|\hat{u}_t||f(x,u_{\theta})-f(x,u)|dx$ can be further integrated over $[0,\tau]\subset[0,T]$ and estimated in a similar manner as in the 6\textsuperscript{th} line of \ref{fEstimation} in Theorem \ref{Q1}. Moreover, using Young's inequality and exploiting the fact that $u,\hat{u}\in C^2(\Omega\times[0,T])$, we arrive at
\begin{equation}\label{f}
    \begin{split}
        &2\int_0^{\tau}\int_{\Omega}|\hat{u}_t||f(x,u_{\theta})-f(x,u)|dxdt\\
        &\leq 2\int_0^{\tau}\int_{\Omega}c2^{r-1}|\hat{u}_t|(|u|^r+\frac{1}{r+1}|\hat{u}|^r)|\hat{u}|dxdt\\
        &\leq 2^rc(\|u\|^r_{C^0}+\frac{1}{r+1}\|\hat{u}\|^r_{C^0})\int_0^{\tau}\int_{\Omega}|\hat{u}_t||\hat{u}|dxdt\\
        &\leq 2^{r-1}c(\|u\|^r_{C^0}+\frac{1}{r+1}\|\hat{u}\|^r_{C^0})\left(\int_0^{\tau}\int_{\Omega}|\hat{u}_t|^2dxdt+\int_0^{\tau}\int_{\Omega}|\hat{u}|^2dxdt\right)\\
        &=:\hat{C}\left(\int_0^{\tau}\int_{\Omega}|\hat{u}_t|^2dxdt+\int_0^{\tau}\int_{\Omega}|\hat{u}|^2dxdt\right).
    \end{split}
\end{equation}
Integrating \ref{Zwischenergebnis1} over $[0,\tau]\subset[0,T]$ and using Cauchy-Schwarz and Young's inequality gives us
\begin{equation*}
    \begin{split}
        \int_{\Omega}|\hat{u}_t(x,\tau)|^2dx\leq&\|\mathcal{R}_{u_1}\|_{L^2(\Omega)}^2+\int_0^{\tau}\int_{\Omega}|\hat{u}_t|^2dx+\|\mathcal{R}_{PDE}\|_{L^2(\Omega_T)}^2+\|\mathcal{R}_{\nabla u}\|^2_{L^2(\Omega)}\\
        &-\int_{\Omega}|\nabla\hat{u}(x,\tau)|^2dx+\hat{C}(\int_0^{\tau}\int_{\Omega}|\hat{u}_t|^2dxdt+\int_0^{\tau}\int_{\Omega}|\hat{u}|^2dxdt)\\
        &+2\sqrt{T|\partial\Omega|}\|\hat{u}\|_{C^1}\|\mathcal{R}_{s,u_t}\|_{L^2(\partial\Omega\times[0,T])}.\\
    \end{split}
\end{equation*}
Rearranging this results leads to 
\begin{equation}\label{Zwischenergebnis2}
    \begin{split}
        &\int_{\Omega}|\hat{u}_t(x,\tau)|^2dx+\int_{\Omega}|\nabla\hat{u}(x,\tau)|^2dx\\
        &\leq\|\mathcal{R}_{t_1}\|_{L^2(\Omega)}^2+\int_0^{\tau}\int_{\Omega}|\hat{u}_t|^2dxdt+\|\mathcal{R}_{PDE}\|_{L^2(\Omega_T)}^2+\|\mathcal{R}_{\nabla u}\|^2_{L^2(\Omega)}\\
        &\quad+\hat{C}(\int_0^{\tau}\int_{\Omega}|\hat{u}_t|^2dxdt+\int_0^{\tau}\int_{\Omega}|\hat{u}|^2dxdt)+2\sqrt{T|\partial\Omega|}\|\hat{u}\|_{C^1}\|\mathcal{R}_{S,u_t}\|_{L^2(\partial\Omega\times[0,T])}.\\
    \end{split}
\end{equation}
Using Corollary \ref{PWICorollary} and Fubini's Theorem we get the following estimate for $\|\hat{u}\|_{L^2(\Omega)}$:
\begin{equation*}
    \|\hat{u}\|_{L^2(\Omega)}\leq C_{pw}\|\nabla\hat{u}\|_{L^2(\Omega)}+\left(\int_{\Omega}|\hat{u}(x,0)|^2dx\right)^{\frac{1}{2}}+\sqrt{T}\left(\int_0^{\tau}\int_{\Omega}|\hat{u}_t|^2dxdt\right)^{\frac{1}{2}},
\end{equation*}
where $C_{pw}$ is the constant we get from the Poincar\'e-Wirtinger inequality.\\
We get
\begin{equation*}
    \|\hat{u}\|^2_{L^2(\Omega)}\leq 2C_{pw}^2\|\nabla\hat{u}\|^2_{L^2(\Omega)}+4\|\mathcal{R}_{u_0}\|^2_{L^2(\Omega)}+4\sqrt{T}\int_0^{\tau}\int_{\Omega}|\hat{u}_t|^2dxdt
\end{equation*}
and
\begin{equation*}
    \frac{1}{2C_{pw}^2}(\|\hat{u}\|^2_{L^2(\Omega)}-4\|\mathcal{R}_{u_0}\|^2_{L^2(\Omega)}-4\sqrt{T}\int_0^{\tau}\int_{\Omega}|\hat{u}_t|^2dxdt)\leq\|\nabla\hat{u}\|^2_{L^2(\Omega)}.
\end{equation*}
Using this inequality in \ref{Zwischenergebnis2} gives us
\begin{equation*}
    \begin{split}
        &\int_{\Omega}|\hat{u}_t(x,\tau)|^2dx+\frac{1}{2C_{pw}^2}\int_{\Omega}|\hat{u}(x,\tau)|^2dx\\
        &\leq\|\mathcal{R}_{u_1}\|_{L^2(\Omega)}^2+\int_0^{\tau}\int_{\Omega}|\hat{u}_t|^2dxdt+\|\mathcal{R}_{PDE}\|_{L^2(\Omega_T)}^2+\|\mathcal{R}_{\nabla u}\|^2_{L^2(\Omega)}+\hat{C}\left(\int_0^{\tau}\int_{\Omega}|\hat{u}_t|^2dxdt+\int_0^{\tau}\int_{\Omega}|\hat{u}|^2dxdt\right)\\
        &\quad+2\sqrt{T|\partial\Omega|}\|\hat{u}\|_{C^1}\|\mathcal{R}_{s,u_t}\|_{L^2(\partial\Omega\times[0,T])}+\frac{2}{C_{pw}^2}\|\mathcal{R}_{u_0}\|_{L^2(\Omega)}^2+\frac{2\sqrt{T}}{C^{pw}_2}\int_0^{\tau}\int_{\Omega}|\hat{u}_t|^2dxdt.
    \end{split}
\end{equation*}
From this it follows that
\begin{equation*}
    \begin{split}
        &\text{min}\{1,\frac{1}{2C_{pw}^2}\}\left(\int_{\Omega}|\hat{u}_t(x,\tau)|^2dx+\int_{\Omega}|\hat{u}(x,\tau)|^2dx\right)\\
        &\leq\|\mathcal{R}_{u_1}\|_{L^2(\Omega)}^2+\|\mathcal{R}_{PDE}\|_{L^2(\Omega_T)}^2+\|\mathcal{R}_{\nabla u}\|^2_{L^2(\Omega)}+2\sqrt{T|\partial\Omega|}\|\hat{u}\|_{C^1}\|\mathcal{R}_{s,u_t}\|_{L^2(\partial\Omega\times[0,T])}\\
        &\quad+\frac{2}{C_{pw}^2}\|\mathcal{R}_{t_0}\|_{L^2(\Omega)}^2+\hat{C}\int_0^{\tau}\int_{\Omega}|\hat{u}|^2dxdt+(1+\hat{C}+\frac{2\sqrt{T}}{C^2_{pw}})\int_0^{\tau}\int_{\Omega}|\hat{u}_t|^2dxdt
    \end{split}
\end{equation*}
and thus
\begin{equation*}
    \begin{split}
        &\int_{\Omega}|\hat{u}_t(x,\tau)|^2dx+\int_{\Omega}|\hat{u}(x,\tau)|^2dx\\
        &\leq\text{max}\{1,2C^2_{pw}\} \biggl( \|\mathcal{R}_{u_1}\|_{L^2(\Omega)}^2+\|\mathcal{R}_{PDE}\|_{L^2(\Omega_T)}^2+\|\mathcal{R}_{\nabla u}\|^2_{L^2(\Omega)}+2\sqrt{T|\partial\Omega|}\|\hat{u}\|_{C^1}\|\mathcal{R}_{s,u_t}\|_{L^2(\partial\Omega\times[0,T])}\\
        &\quad+\frac{2}{C_{pw}^2}\|\mathcal{R}_{u_0}\|_{L^2(\Omega)}^2+(1+\hat{C}+\frac{2\sqrt{T}}{C^2_{pw}})(\int_0^{\tau}\int_{\Omega}|\hat{u}|^2dxdt+\int_0^{\tau}\int_{\Omega}|\hat{u}_t|^2)dxdt \biggr).
    \end{split}
\end{equation*}
Using Gr\"onwall's inequality then gives us
\begin{equation*}
    \begin{split}
        &\int_{\Omega}|\hat{u}_t(x,\tau)|^2dx+\int_{\Omega}|\hat{u}(x,\tau)|^2dx\\
        &\leq\text{max}\{1,2C^2_{pw}\}\biggl(\|\mathcal{R}_{u_1}\|_{L^2(\Omega)}^2+\|\mathcal{R}_{PDE}\|_{L^2(\Omega_T)}^2+\|\mathcal{R}_{\nabla u}\|^2_{L^2(\Omega)}+2\sqrt{T|\partial\Omega|}\|\hat{u}\|_{C^1}\|\mathcal{R}_{s,u_t}\|_{L^2(\partial\Omega\times[0,T])}\\
        &\quad+\frac{2}{C_{pw}^2}\|\mathcal{R}_{u_0}\|_{L^2(\Omega)}^2\biggr)\text{exp}\biggl(T\bigl(\text{max}\{1,2C^2_{pw}\}(1+\hat{C}+\frac{2\sqrt{T}}{C^2_{pw}})\bigr)\biggr).
    \end{split}
\end{equation*}
Integrating over $[0,T]$ leads to the result
\begin{equation*}
    \begin{split}
        &\int_0^T\int_{\Omega}|\hat{u}_t(x,t)|^2dxdt+\int_0^T\int_{\Omega}|\hat{u}(x,t)|^2dxdt\\
        &\leq T\text{max}\{1,2C^2_{pw}\}\biggl(\|\mathcal{R}_{u_1}\|_{L^2(\Omega)}^2+\|\mathcal{R}_{PDE}\|_{L^2(\Omega_T)}^2+\|\mathcal{R}_{\nabla u}\|^2_{L^2(\Omega)}+2\sqrt{T|\partial\Omega|}\|\hat{u}\|_{C^1}\|\mathcal{R}_{s,u_t}\|_{L^2(\partial\Omega\times[0,T])}\\
        &\quad+\frac{2}{C_{pw}^2}\|\mathcal{R}_{u_0}\|_{L^2(\Omega)}^2\biggr)\text{exp}\biggl(T\bigl(\text{max}\{1,2C^2_{pw}\}(1+\hat{C}+\frac{2\sqrt{T}}{C^2_{pw}})\bigr)\biggr)
    \end{split}.
\end{equation*}
\end{proof}

\begin{proof}[Proof of Theorem \ref{Q3}]
    From the quadrature rule \ref{Quadrature2} it follows that for $i\in\{u_0,u_1,\nabla u\}$
    \begin{equation*}
    \begin{split}
        &\|\mathcal{R}_{PDE}\|^2_{L^2(\Omega_T)}\leq C_{\Omega_T}\|\mathcal{R}_{PDE}^2\|_{C^2}M^{-\frac{2}{d+1}}+\mathcal{Q}^{PDE}_{M_{PDE}}[\mathcal{R}_{PDE}^2],\\
        &\|\mathcal{R}_{s,u_t}\|_{L^2(\partial\Omega\times[0,T])}\leq\sqrt{ C_{\partial\Omega\times[0,T]}\|\mathcal{R}_{s,u_t}^2\|_{C^2}}M^{-\frac{1}{d}}+\sqrt{\mathcal{Q}^{s,u_t}_{M_t}[\mathcal{R}_{s,u_t}^2]},\\
        &\|\mathcal{R}_i\|^2_{L^2(\Omega)}\leq C_{\Omega}\|\mathcal{R}_i^2\|_{C^2}M^{-\frac{2}{d}}+\mathcal{Q}^i_{M_t}[\mathcal{R}_i^2].
    \end{split}
    \end{equation*}
Combining this result with Theorem \ref{Q2} gives us the main result of Theorem \ref{Q3}. Additionally, we note that Lemma \ref{CnNorm} gives us the following bound on the $C^n$-norm of $u_{\theta}$
\begin{equation*}
    \|u_{\theta}\|_{C^n}\leq 16^L(d+1)^{2n}\bigl(e^2n^4W^3R^n\|\sigma\|_{C^n}\bigr)^{nL},
\end{equation*}
which we can use to bound $\|u_{\theta}-u\|_{C^1}\leq\|u_{\theta}\|_{C^1}+\|u\|_{C^1}$ as well as deriving the following constants:
\begin{equation}\label{ConstantsMidPointRule}
    \begin{split}
        C_1&:=C_{\Omega}\|\mathcal{R}_{u_1}^2\|_{C^2}\leq 4C_{\Omega}\|\mathcal{R}_{u_1}\|_{C^2}^2\leq 4C_{\Omega}(\|u_{\theta_t}\|_{C^2}+\|u_t\|_{C^2})^2\leq 8C_{\Omega}(\|u_{\theta}\|_{C^3}^2+\|u\|_{C^3}^2)\\
        &\leq 8C_{\Omega}\bigl(16^{2L}(d+1)^{12}(e^23^4W^3R^3\|\sigma\|_{C^3})^{6L}+\|u\|_{C^3}^2\bigr),\\
        \\[-5mm]
        C_2&:=C_{\Omega_T}\|\mathcal{R}_{PDE}^2\|_{C^2}\leq 4C_{\Omega_T}\|\mathcal{R}_{PDE}\|^2_{C^2}\leq 4C_{\Omega_T}\left(\|u_{\theta_{tt}}\|_{C^2}+\|\Delta u_{\theta}\|_{C^2}+\|au_{\theta_t}\|_{C^2}+\|f(x,u_{\theta})\|_{C^2}\right)^2\\
        &\leq 16 C_{\Omega_T}\left(\|u_{\theta_{tt}}\|_{C^2}^2+\|\Delta u_{\theta}\|_{C^2}^2+\|au_{\theta_t}\|_{C^2}^2+\|f(x,u_{\theta})\|_{C^2}^2\right)\\
        &\leq 16 C_{\Omega_T}\left(\|u_{\theta}\|_{C^4}^2+d\| u_{\theta}\|_{C^4}^2+16\|a\|_{C^2}^2\|u_{\theta}\|_{C^3}^2+c^2\mathrm{max}\{\|u_{\theta}\|_{C^0}^{2(r+1)},\|u_{\theta}\|_{C^0}^{2r}, \|u_{\theta}\|_{C^0}^{2(r-1)}\}\right)\\
        &\leq 16 C_{\Omega_T}\biggl(16^{2L}(d+1)^{17}\bigl(e^24^4W^3R^4\|\sigma\|_{C^4}\bigr)^{8L}+16\|a\|_{C^2}^216^{2L}(d+1)^{12}\bigl(e^23^4W^3R^3\|\sigma\|_{C^3}\bigr)^{6L}\\
        &\quad+c^22^rR^{r+1}(W^{r+1}\|\sigma\|^{r+1}_{C^0}+1)\biggr),\\
        \\[-5mm]
        C_3&:=C_{\Omega}\|\mathcal{R}_{\nabla u}^2\|_{C^2}\leq 4C_{\Omega}\|\mathcal{R}_{\nabla u}\|_{C^2}^2\leq 4C_{\Omega}(\|\nabla u_{\theta}\|_{C^2}+\|\nabla u\|_{C^2})^2\leq 8C_{\Omega}(\|\nabla u_{\theta}\|_{C^2}^2+\|\nabla u\|_{C^2}^2)\\
        &\leq 8dC_{\Omega}(\|u_{\theta}\|_{C^3}^2+\|u\|_{C^3}^2)\leq8dC_{\Omega}\bigl(16^{2L}(d+1)^{12}(e^23^4W^3R^3\|\sigma\|_{C^3})^{6L})+\|u\|_{C^3}^2\bigr),\\
        \\[-5mm]
        C_4&:=\sqrt{C_{\partial\Omega\times[0,T]}\|\mathcal{R}_{s,u_t}^2\|_{C^2}}\leq2\sqrt{C_{\partial\Omega\times[0,T]}}\|\mathcal{R}_{s,u_t}\|_{C^2}=2\sqrt{C_{\partial\Omega\times[0,T]}}\|u_{\theta_t}\|_{C^2}\leq2\sqrt{C_{\partial\Omega\times[0,T]}}\|u_{\theta}\|_{C^3}\\
        &\leq2\sqrt{C_{\partial\Omega\times[0,T]}}16^L(d+1)^6(e^23^4W^3R^3\|\sigma\|_{C^3})^{3L},\\
        \\[-5mm]
        C_5&:=\frac{2C_{\Omega}}{C^2_{pw}}\|\mathcal{R}_{u_0}^2\|_{C^2}\leq\frac{8C_{\Omega}}{C^2_{pw}}\|\mathcal{R}_{u_0}\|^2_{C^2}\leq\frac{8C_{\Omega}}{C^2_{pw}}(\|u_{\theta}\|_{C^2}+\|u\|_{C^2})^2\leq\frac{16C_{\Omega}}{C^2_{pw}}(\|u_{\theta}\|^2_{C^2}+\|u\|^2_{C^2})\\
        &\leq\frac{16C_{\Omega}}{C^2_{pw}}\bigl(16^{2L}(d+1)^8(e^22^4W^3R^2\|\sigma\|_{C^2})^{4L}+\|u\|^2_{C^2}),
    \end{split}
\end{equation}
where we used that for every residual it holds that $\|\mathcal{R}_i^2\|_{C^n}\leq 2^n\|\mathcal{R}_i\|^2_{C^n}$ by Lemma \ref{A6}. Moreover, for the bound on $C_2$ we used Lemma \ref{A6}, Lemma \ref{C0Norm} and the the fact that the derivatives of $f$ are bounded according to assumption (A2).
\end{proof}

\begin{proof}[Proof of Theorem \ref{apriori}]
    Let $L=3$ and $N\in\mathbb{N}, N>5$. From Theorem \ref{Q1} and Corollary \ref{BoundGeneralizationError} we know that there exists $\hat{\theta}\in\Theta$ such that $\mathcal{E}_G(\hat{\theta})=\mathcal{O}(\mathrm{ln^2}(N)N^{-k+1})$ and $W=CN^{d+1}$ (the exact value of $C$ follows from theorem \ref{Q1}), plus $R=\mathcal{O}(N\mathrm{ln}(N))$ for $n$ large enough.\\
    Let $\theta^*(\mathcal{S})$ be the solution to \ref{OptimizationProblem2}. Since in particular $k\geq4$ and $u\in H^{k+1}(\Omega\times(0,T))$ the Sobolev Embedding Theorem, see, e.g. \cite{Brezis}, implies that $u\in C^3(\Omega\times(0,T))$ and we can use Theorem \ref{Q3} and the assumption $24L-1\geq r$ to derive
    \begin{equation*}
        \|u-u_{\theta^*(\mathcal{S})}\|^2_{L^2(\Omega\times[0,T])}+\|u_t-u_{\theta^*(\mathcal{S})_t}\|^2_{L^2(\Omega\times[0,T])}\lesssim (W^3R^4)^{8L}\left(\mathcal{E}_T(\theta^*(\mathcal{S}),\mathcal{S})^2+M_t^{-\frac{2}{d}}+M_{PDE}^{-\frac{2}{d+1}}+M_s^{-\frac{1}{d}}\right).
    \end{equation*}
    Since $\theta^*(\mathcal{S})$ is the global minimum, it holds that $\mathcal{E}_T(\theta^*(\mathcal{S}),\mathcal{S})\leq \mathcal{E}_T(\hat{\theta}(\mathcal{S}),\mathcal{S})$. Using Corollary \ref{BountTraininError} we can bound $\mathcal{E}_T(\hat{\theta}(\mathcal{S}),\mathcal{S})$ by $\mathcal{E}_G(\hat{\theta})$ and get
    \begin{equation*}
        \|u-u_{\theta^*(\mathcal{S})}\|^2_{L^2(\Omega\times[0,T])}+\|u_t-u_{\theta^*(\mathcal{S})_t}\|^2_{L^2(\Omega\times[0,T])}\lesssim (W^3R^4)^{8L}\left(\mathcal{E}_G(\hat{\theta})^2+M_t^{-\frac{2}{d}}+M_{PDE}^{-\frac{2}{d+1}}+M_s^{-\frac{1}{d}}\right).
    \end{equation*}
    We can rewrite everything in terms of $N$ and use $\mathrm{ln}^2(N)<N$ for $N>5$ to arrive at
    \begin{equation*}
        \begin{split}
            &\|u-u_{\theta^*(\mathcal{S})}\|^2_{L^2(\Omega\times[0,T]}+\|u_t-u_{\theta^*(\mathcal{S})_t}\|^2_{L^2(\Omega\times[0,T]}\\
            &\lesssim\left(N^{3(d+1)}\mathrm{ln}^4(N)N^4\right)^{24}\left(\mathrm{ln}^4(N)N^{2(-k+1)}+M_t^{-\frac{2}{d}}+M_{PDE}^{-\frac{2}{d+1}}+M_s^{-\frac{1}{d}}\right)\\
            &\lesssim N^{24(3d+9)}\left(N^{-2k+4}+M_t^{-\frac{2}{d}}+M_{PDE}^{-\frac{2}{d+1}}+M_s^{-\frac{1}{d}}\right).
        \end{split}
    \end{equation*}
    We choose $N$ and the training set sizes $M_t, M_{PDE}$ and $M_s$ such that 
    \begin{equation*}
        N^{24(3d+9)}\left(N^{-2k+4}+M_t^{-\frac{2}{d}}+M_{PDE}^{-\frac{2}{d+1}}+M_s^{-\frac{1}{d}}\right)\overset{!}{=}\mathcal{O}(\varepsilon^2),
    \end{equation*}
    which leads to $N=\varepsilon^{-\frac{1}{k-\eta}}, M_{PDE}=\varepsilon^{-\frac{-(d+1)(\eta-1)}{k-\eta}}, M_t=\varepsilon^{-\frac{-d(\eta-1)}{k-\eta}}$ and $M_s=\varepsilon^{-\frac{-2d(\eta-1)}{k-\eta}}$. For example, for the first term we want
    \begin{equation*}
        N^{24(3d+9)-2k+4}=\varepsilon^2,
    \end{equation*}
    which leads to $N=\varepsilon^{\frac{1}{2(18d+55)-k}}=:\varepsilon^{\frac{1}{\eta-k}}$ and therefore $\varepsilon=N^{\eta-k}$. For $\varepsilon\ll 1$ this leads to the condition $k>\eta$.
\end{proof}

\section*{Acknowledgments}
This work of G. Kutyniok was supported in part by the ONE Munich Strategy Forum (LMU Munich, TU Munich, and the Bavarian Ministery for Science and Art).

G. Kutyniok acknowledges support from the Konrad Zuse School of Excellence in Reliable AI (DAAD), the Munich Center for Machine Learning (BMBF) as well as the German Research Foundation under Grants DFG-SPP-2298, KU 1446/31-1 and KU 1446/32-1 and under Grant DFG-SFB/TR 109, Project C09 and the Federal Ministry of Education and Research under Grant MaGriDo.

 \bibliographystyle{elsarticle-num} 
 \bibliography{references}





\end{document}